\documentclass[10pt]{amsart}

\usepackage{amssymb,array}

\newtheorem{thm}{Theorem}

\newtheorem{lem}[thm]{Lemma}
\newtheorem{cor}[thm]{Corollary}

\newtheorem{prop}[thm]{Proposition}

   \newtheorem{lemma}[thm]{Lemma}  
\newtheorem{theorem}[thm]{Theorem}  

\theoremstyle{definition}
\newtheorem{defn}[thm]{Definition}

\newtheorem{say}[thm]{}
\newtheorem{exmp}[thm]{Example}

\newtheorem{rem}[thm]{Remark}

\newtheorem{defn-thm}[thm]{Definition--Theorem}  
\newtheorem{defn-lem}[thm]{Definition--Lemma}  

\theoremstyle{remark}

\renewcommand{\c}[0]{{\mathbb C}}  

\renewcommand{\o}[0]{{\mathcal O}} 

\newcommand{\q}[0]{{\mathbb Q}}
\newcommand{\map}[0]{\dasharrow}
\newcommand{\qtq}[1]{\quad\mbox{#1}\quad}

\newcommand{\pic}[0]{\operatorname{Pic}}
\newcommand{\gal}[0]{\operatorname{Gal}}

\newcommand{\Hom}[0]{\operatorname{Hom}}

\newcommand{\aut}[0]{\operatorname{Aut}}

\newcommand{\sing}[0]{\operatorname{Sing}}    
\newcommand{\ex}[0]{\operatorname{Ex}}

\newcommand{\isom}[0]{\operatorname{Isom}}
\newcommand{\alb}[0]{\operatorname{alb}}
\newcommand{\Alb}[0]{\operatorname{Alb}}

\newcommand{\stab}[0]{\operatorname{Stab}}
\newcommand{\simq}[0]{\sim_{\q}}
\newcommand{\age}[0]{\operatorname{age}}

\renewcommand{\bar}{\overline}

\newcommand{\C}{{\mathbb{C}}}

\newcommand{\Q}{{\mathbb{Q}}}
\newcommand{\Z}{{\mathbb{Z}}}

\newcommand{\End}{\mathrm{End}}
\newcommand{\Aut}{\mathrm{Aut}}

\newcommand{\Gal}{\mathrm{Gal}}
\newcommand{\GL}{\mathrm{GL}}
\newcommand{\PGL}{\mathrm{PGL}}

\newcommand{\Span}{\mathrm{Span}}

\newcommand{\tr}{\mathrm{tr}}

\makeindex

\begin{document}
\bibliographystyle{amsalpha}

\title{Quotients of Calabi-Yau varieties}\author{J\'anos Koll\'ar}
\email{kollar@math.princeton.edu}
\address{Department of Mathematics \\
	Princeton University \\
	Fine Hall, Washington Road Princeton, NJ 08544-1000 \\
	U.S.A.}

\author{Michael Larsen}
\email{larsen@math.indiana.edu}
\address{Department of Mathematics\\
    Indiana University \\
    Bloomington, IN 47405\\
    U.S.A.}

\date{January 20, 2007}

\maketitle

Let $X$ be a Calabi-Yau variety over $\C$, that is, a  projective variety
with canonical singularities 
whose canonical class is numericaly
trivial. Let $G$ be a finite group acting on $X$ and consider the
quotient variety $X/G$. The aim of this paper is to determine the place of
$X/G$ in the birational classification of varieties. That is,
we determine the Kodaira dimension of $X/G$ and decide when
it is uniruled or rationally connected.

If $G$ acts without fixed points, then
$\kappa(X/G)=\kappa(X)=0$, thus the interesting case is
when $G$ has fixed points. 
We answer the above questions in terms of the
action of the stabilizer subgroups near the fixed points.

The answer is especially nice if $X$ is smooth. In the introduction we 
concentrate on this case. The precise general results are
formulated later.

\begin{defn}\label{rt.defn}
Let $V$ be a complex vector space and $g\in \GL(V)$ an
element of finite order. Its eigenvalues (with multiplicity) can be 
written as $e(r_1),\ldots,e(r_n)$,
where $e(x) := e^{2\pi i x}$ and $0\le r_i<1$.
Following \cite{it-re, reid-mckay}, we define the {\it age} of $g$ as
$$
\age(g):=r_1+\cdots+r_n.
$$
Let $G$ be a finite group and $(\rho,V)$ a finite dimensional 
complex representation of $G$.
We say that $\rho:G\to \GL(V)$  satisfies 
the (local) {\it Reid-Tai condition} if
$\age(\rho(g))\geq 1$ for every $g\in G$ for which $\rho(g)$ is not the identity
(cf.\ \cite[3.1]{reid}).

Let $G$ be a finite group acting on  a smooth 
 projective variety $X$. We say that
the $G$-action satisfies the (global) {\it Reid-Tai condition}
if for every $x\in X$, the stablizer representation
$\stab_x(G)\to \GL(T_xX)$ satisfies 
the (local) Reid-Tai condition.
\end{defn}

Our first result relates the uniruledness of $X/G$
to the Reid-Tai condition.

\begin{thm}\label{intr.quot.thm} Let $X$ be a smooth 
 projective Calabi-Yau variety
and $G$ a finite group acting on $X$. The following are equivalent:
\begin{enumerate}
\item $\kappa(X/G)=0$.
\item $X/G$ is not uniruled.
\item The $G$-action satisfies the global Reid-Tai condition.
\end{enumerate}
\end{thm}

It is conjectured in general that being uniruled is equivalent to
having  Kodaira dimension  $-\infty$. As  part of Theorem 2, we establish this
equivalence for varieties of the form $X/G$.

It can happen that $X/G$ is uniruled but not rationally
connected. The simplest example
is when $X=X_1\times X_2$ is a product, $G$ acts trivially
on $X_1$ and $X_2/G$ is rationally connected.
  Then $X/G\cong X_1\times (X_2/G)$
is a product of the Calabi--Yau variety $X_1$ and of the
rationally connected variety $X_2/G$.
We show that this is essentially the only way
that $X/G$ can be uniruled but not rationally connected.
The key step is a  description of  rational maps
from Calabi-Yau varieties to lower dimensional nonuniruled varieties.

\begin{thm}\label{intr.map.thm} Let $X$ be a smooth, simply connected 
 projective Calabi-Yau variety 
and $f:X\map Y$ a dominant map such that $Y$ is not uniruled.
Then one can write
$$
f: X\stackrel{\pi}{\longrightarrow} X_1\stackrel{g}{\map} Y
$$
where $\pi$ is a projection to a direct factor of $X\cong X_1\times X_2$
and $g:X_1\map Y$ is generically finite.
\end{thm}

Applying this to the MRC-fibration of $X/G$,
 we obtain the following.

\begin{cor}\label{intr.quot.thm2} Let $X$ be a smooth, simply connected 
 projective Calabi-Yau variety which is not a nontrivial product of
two Calabi-Yau varieties.
Let $G$ be a finite group acting on $X$. The following are equivalent:
\begin{enumerate}
\item $X/G$ is uniruled.
\item $X/G$ is rationally connected.
\item The $G$-action does not satisfy the global Reid-Tai condition.
\end{enumerate}
\end{cor}

Next we turn our attention to a study of the local
Reid-Tai condition. For any given representation it is
relatively easy to decide if the  Reid-Tai condition is satisfied or not.
It is, however, quite difficult to get a good understanding
of all representations that satisfy it. For instance, it is quite
tricky to determine all $\leq 4$-dimensional  representations of cyclic groups
that satisfy the Reid-Tai condition, 
cf.\ \cite{mo-st, mor-q, momomo, reid-yp}.
These turn out to be rather special.

By contrast, we claim that every representation of a ``typical''
nonabelian group satisfies the Reid-Tai condition.
The groups which have some representation violating the Reid-Tai
condition are closely related to complex reflection groups.
In the second part of the paper we provide a kit for building all 
of them, using basic building blocks, all but finitely
many of which are (up to projective equivalence) reflection groups.

Let $G$ be a finite group and $(\rho,V)$ a finite dimensional 
complex representation of $G$ such that 
 $(\rho,V)$ does not  satisfy the  (local) Reid-Tai condition.
That is,  $G$ has an element $g$ such that $0<\age(\rho(g))<1$.
We say that such a pair $(G,V)$ is a {\it non-RT pair} and $g$ is an \emph{exceptional element}.
There is no essential gain in generality in allowing 
$\rho\colon G\to \GL(V)$ not to be faithful.
We therefore assume that $\rho$ is faithful, and remove it from the notation,
regarding
$G$ as a subgroup of $\GL(V)$ (which is to be classified up to conjugation).
If the conjugacy class of $g$ does not generate the full group $G$, 
it must generate a normal subgroup $H$ of $G$ such that $(H,V)$ is again a non-RT pair.
After classifying the cases for which the conjugacy class of $g$ 
generates $G$, we can
take the normalization of each such $G$ in $\GL(V)$; all finite subgroups 
intermediate between $G$
and this normalizer give further examples.   
If $V$ is reducible, then for every irreducible factor $V_i$ of 
$V$ on which $g$ acts non-trivially, $(G,V_i)$ is again a non-RT pair with exceptional
element $g$.  Moreover, if the conjugacy class of $g$ generates $G$, then $g$
must be an exceptional element for every non-trivial factor $V_i$ of $V$.
These reduction steps motivate the following definition:

\begin{defn}
A {\it basic non-RT pair} is an ordered pair $(G,V)$ consisting of a 
finite group $G$ and a faithful irreducible representation $V$ such that the conjugacy 
class of any 
exceptional element $g\in G$ generates $G$.
\end{defn}

Given a basic non-RT pair $(G,V)$ and a positive integer $n$, we define
$G_n = G\times \Z/n\Z$ and let $V_n$ denote the tensor product of $V$ with the
character of $\Z/n\Z$ sending $1$ to $e(1/n)$.  Then $V_n$ is always an irreducible representation
of $G_n$ and is faithful if $n$ is prime to the order of $G$.  Also if the conjugacy class of $g\in G$
generates the whole group, $a\in\Z$ is relatively prime to $n$, and $n$ is prime to the order of $G$,
then the conjugacy class of $(g,\bar a)\in G_n$ generates $G_n$.  Finally, if
$g$ is an exceptional element, $(a,n)=1$, and $a/n>0$ is sufficiently small,
then $(g,\bar a)$ is exceptional.  Thus for each
basic non-RT pair $(G,V)$, there are infinitely many other basic non-RT pairs which are projectively
equivalent to it.  To avoid this complication, we seek to classify basic non-RT pairs only up to 
projective equivalence:

\begin{defn}
Pairs $(G_1,V_1)$ and $(G_2,V_2)$
are \emph{projectively equivalent} if there exists 
an isomorphism $\PGL(V_1)\to\PGL(V_2)$ mapping the image
of $G_1$ in $\PGL(V_1)$ isomorphically to the image of $G_2$ in $\PGL(V_2)$.
\end{defn}

We recall that a \emph{pseudoreflection} $g\in\GL(V)$ is an element of finite order which
fixes a codimension $1$ subspace of $V$ pointwise.  A (complex) reflection group
is a finite subgroup of $\GL(V)$ which is generated by pseudoreflections.  We say that
$(G,V)$ is a reflection group if $G$ is a reflection group in $\GL(V)$.

\begin{defn}
We say that a basic non-RT pair $(G,V)$ is of {\it reflection type}
if $(G,V)$ is projectively equivalent to some reflection group $(G',V')$.  
\end{defn}

The reflection groups are classified in \cite{st}.
Note that every pseudoreflection is of exceptional type, so 
every reflection group is a non-RT pair.  On the other hand, there
may be elements of exceptional type in a reflection group which are \emph{not} pseudoreflections.
Moreover, not every irreducible reflection group gives rise to a basic non-RT pair; it 
may happen that a particular conjugacy class of pseudoreflections fails to generate
the whole group, since some irreducible reflection groups have multiple conjugacy classes
of pseudoreflections.   

We are interested in the other direction, however, and here we have the following theorem.

\begin{theorem}\label{nonRT.almost.class}
Up to projective equivalence, there are only finitely 
many basic non-RT pairs which are \emph{not} of reflection type.
\end{theorem}

We know of only one example.   In principle, our proof provides an effective way 
(via the classification of finite simple groups) to determine all examples.

Finally, we study in detail the case when $X=A$ is an Abelian variety.
In fact, this case was the starting point of our investigations.
For Abelian varieties, the induced representations
$\stab_xG\to V=T_xA$ have the property that
$V+V^*$ is isomorphic to the rational representation
of $\stab_xG$ on 
 $H^1(A,\q)$. This property substantially reduces the number of cases that we need
to consider and allows us to show that a 
basic non-RT pair arising from an Abelian variety which is not of reflection type is of the
unique known type.
A more precise statement is given in Theorem~\ref{AV.case} below.
Unlike the proof of Theorem~\ref{nonRT.almost.class}, the proof of this theorem
does not make use of the classification of finite simple groups.

A consequence of the analysis which ultimately gives Theorem~\ref{AV.case} is the following easier statement.

\begin{thm} 
\label{at.least.four}
Let $A$ be a simple Abelian variety of dimension
$\geq 4$ and $G$ a finite group acting on $A$. Then
$A/G$ has canonical singularities and $\kappa(A/G)=0$.
\end{thm}

\noindent\textbf{Acknowledgments.} We thank Ch.\ Hacon and  S.\ Kov\'acs
for hepful comments and references
and De-Qi Zhang for correcting an error in an earlier  version
of Theorem \ref{intr.map.thm}. 
Partial financial support to JK and ML was provided by  the NSF under grant numbers
DMS-0500198 and DMS-0354772. 

\section{Uniruled quotients}

Let $X$ be a  projective Calabi--Yau variety
and $G\subset \aut(X)$
a finite group of automorphism.
Our aim is to decide when is the quotient variety $X/G$
uniruled or rationally connected.
Our primary interest is in the case when $X$ is smooth,
but the proof works without change whenever $X$ has canonical singularities
and $K_X$ is numerically trivial.
(Note that by \cite[8.2]{kaw}, $K_X$  numerically trivial implies that
$K_X$ is torsion.)

\begin{thm}\label{cy.quot.ur.thm} Let $X$ be a  projective Calabi--Yau variety
and $G$ a finite group acting on $X$. The following are equivalent:
\begin{enumerate}
\item  $G$ acts freely outside a codimension $\geq 2$ set and
$X/G$ has canonical singularities.
\item $\kappa(X/G)=0$.
\item $X/G$ is not uniruled.
\end{enumerate}
\end{thm}

Proof. Set $Z:=X/G$ and let $D_i\subset Z$ be the branch divisors
of the quotient map $\pi:X\to Z$ with branching index $e_i$.
Set $\Delta:=\sum \bigl(1-\frac1{e_i}\bigr)D_i$.
By the Hurwitz formula,
$$
K_X\simq \pi^*\bigl(K_Z+\Delta\bigr),
$$
where $\simq$ means that some nonzero integral  multiples of the two sides are
linearly equivalent.
Thus 
$K_Z+\Delta\simq 0$ and hence 
Theorem~\ref{cy.quot.ur.thm}) is a special case of Theorem~\ref{logcy.ur.thm}.

\begin{thm}\label{logcy.ur.thm} Let $Z$ be a  projective variety and
$\Delta$ an effective divisor on $Z$ such that
$K_Z+\Delta\simq 0$. The following are equivalent:
\begin{enumerate}
\item  $\Delta=0$  and
$Z$ has canonical singularities.
\item $\kappa(Z)=0$.
\item $Z$ is not uniruled.
\end{enumerate}
\end{thm}

Proof. If $\Delta=0$ then $K_Z$ is numerically trivial.
Let $g:Y\to Z$ be a resolution of singularities
and write
$$
K_Y\simq g^*K_Z+\sum a_iE_i\simq \sum a_iE_i
$$
where the $E_i$ are $g$-exceptional and $a_i\geq 0$ for every $i$ iff
$Z$ has canonical singularities.

Thus if (1) holds then $K_Y\simq \sum a_iE_i$ is effective
and so $\kappa(Y)\geq 0$. Since $\sum a_iE_i$ is exceptional,
no multiple of it moves, hence $\kappa(Y)= \kappa(Z)=0$.

The implication (2) $\Rightarrow$ (3) always holds
(cf.\ \cite[IV.1.11]{rc-book}). It is conjectured that
in fact (2) is equivalent to (3), but this is known only in
dimensions $\leq 3$ (cf.\ \cite[3.12--13]{km-book}).

Thus it remains to prove that if (1) fails then $Z$ is
uniruled. We want to use the Miyaoka-Mori criterion \cite{mi-mo}
to get uniruledness.
That is, a  projective variety $Y$ is uniruled if
an open subset of 
it is covered by  projective   curves $C\subset Y$  such that
$K_Y\cdot C<0$ and $C\cap \sing Y=\emptyset$.

If $\Delta\neq 0$
then we can take $C\subset Z$ to be any smooth complete intersection
curve which does not intersect the singular locus of $Z$.

Thus we can assume from now on that $\Delta=0$
and, as noted before, 
$K_Y\simq \sum a_iE_i$
where the $E_i$ are $g$-exceptional and $a_i< 0$ for some $i$ 
since
$Z$ does not have canonical singularities
by assumption.
For notational convenience assume that
 $a_1<0$.  

Ideally, we would like  to find curves 
$C\subset Y$ such that $C$ intersects $E_1$
but no other $E_i$.  If such a $C$ exists then
$$
(K_Y\cdot C)=(\sum a_iE_i\cdot C) =a_1(E_1\cdot C)<0.
$$
We are not sure that such curves exist.
(The condition $K_Z\equiv 0$ puts strong restrictions
on the singularities of $Z$ and creates a rather special
situation.)

Fortunately, it is sufficient 
to find curves $C$ such that $(C\cdot E_1)$
is big and the other $(C\cdot E_i)$ are small. This is enough to give
$K_Y\cdot C<0$.

\begin{lem} \label{ex.nonef.nonpef} Let $Z$ be a normal projective variety
over a field of arbitrary characteristic, 
 $g:Y\to Z$  a birational morphism
and $E=\sum a_iE_i$ a noneffective $g$-exceptional Cartier divisor.
Then $Y$ is covered by curves $C$ such that $(E\cdot C)<0$.
\end{lem}

Proof.  Our aim is to reduce the problem to a carefully
chosen surface $S\to Y$ and then construct such curves
directly on $S$.  At each step we make sure that $S$ can be chosen
to pass through any general point of $Y$, so if we can cover
these surfaces $S$ with curves $C$, the resulting curves also cover $Y$.

Assume that $a_1<0$.  If $\dim g(E_1)>0$ then we can cut $Y$
by pull backs of hypersurface sections of $Z$
 and use induction on the dimension.

Thus assume that $g(E_1)$ is a point. The next step
would be to cut with hypersurface sections of $Y$.
 The problem is that in this process
some of the divisors $E_i$ may become nonexceptional, even ample.
Thus first we need to kill all the other $E_i$ such that
$\dim g(E_i)>0$.

To this end, construct a series of varieties $\sigma_i:Z_i\to Z$
starting with $\sigma_0:Z_0= Z$ as follows.
Let $W_i\subset Z_i$ be the closure of
$(\sigma_i^{-1}\circ g)_*(E_1)$, $\pi_i:Z_{i+1}\to Z_i$ 
 the blow up of $W_i\subset Z_i$ and $\sigma_{i+1}=\sigma_i\circ \pi_i$.

By Abhyankar's lemma (in the form given in \cite[2.45]{km-book}),
 there is an index  $j$
such that if $Y'\subset Z_j$ denotes the main component
and   $g':=\sigma_j:Y'\to Z$ the induced 
 birational morphism then  the following hold:
\begin{enumerate}
\item $g'$ is an isomorphism over $Z\setminus g(E_1)$, and
\item $h:=g^{-1}\circ g':Y'\map Y$ is a local isomorphism
over the generic point of $E_1$. 
\end{enumerate}
Thus $h^*(\sum a_iE_i)$ is $g'$-exceptional and
not effective (though it is only guaranteed to be Cartier
outside the indeterminacy locus of $h$).

Now we can cut by hypersurface sections of $Y'$
to get a surface $S'\subset Y'$. Let
$\pi:S\to S'$ be a resolution 
such that $h\circ \pi:S\to Y$ is a morphism
and $f:=(g'\circ \pi):S\to T:=g'(S')\subset Z$
the induced morphism.
Then $(h\circ \pi)^*E$ is exceptional over $T$
and not effective.
Thus it is sufficient to  prove Lemma~\ref{ex.nonef.nonpef} in case
$S=Y$ is a smooth surface.

Fix an ample divisor $H$ on $S$ such that
$$
H^1(S, \o_S(K_S+H+L))=0
$$
for every nef divisor $L$. (In characteristic 0
any ample divisor works by the  Kodaira vanishing theorem.
In positive characteristic, one can use for instance
\cite[Sec.II.6]{rc-book}  to show that any $H$ such that
$(p-1)H-K_S-4(\mbox{some ample divisor})$ is nef
has this property.)

Assume next that $H$ is also very ample and pick $B\in |H|$.
Using the 
exact sequence  
$$
0\to \o_S(K_S+2H+L)\to \o_S(K_S+3H+L)\to \o_B(K_B+2H|_B+L|_B)\to 0,
$$
we conclude that $\o_S(K_S+3H+L)$ is generated by global sections.

By the Hodge index theorem, 
the intersection product on the curves $E_i$ is negative definite
hence nondegenerate.
Thus we can find a linear combination $F=\sum b_iE_i$ such that
$F\cdot E_1>0$ and $F\cdot E_i=0$ for every $i\neq 1$.
Choose $H_Z$ ample on $Z$ such that $F+g^*H_Z$ is nef.

Thus the linear system $|K_S+3H+m(F+g^*H_Z)|$ is base point free
for every $m\geq 0$. Let $C_m\in |K_S+3H+m(F+g^*H_Z)|$
be a general irreducible curve. Then
$$
\begin{array}{rcl}
(C_m\cdot E_1)&=&m(F\cdot E_1)+(\mbox{constant}),\qtq{and}\\
(C_m\cdot E_i)&=&(\mbox{constant})\qtq{for $i>1$.}
\end{array}
$$
Thus $(C_m\cdot E)\to -\infty$ as $m\to \infty$.\qed
\medskip

The following consequence of Lemma~\ref{ex.nonef.nonpef}
is of independent interest. 
In characteristic zero, Corollary~\ref{ex.nonef.nonpef.cor}
is equivalent to  Lemma~\ref{ex.nonef.nonpef} by  \cite{d-etal},
thus one can use Corollary~\ref{ex.nonef.nonpef.cor} to give an
alternate proof of Lemma~\ref{ex.nonef.nonpef}.
The equivalence should also hold in any characteristic.

\begin{cor}[Lazarsfeld, unpublished] \label{ex.nonef.nonpef.cor} 
Let $Z$ be a normal, projective variety, 
 $g:Y\to Z$  a birational morphism
and $E=\sum a_iE_i$ a $g$-exceptional  Cartier divisor.
Then $E $ is pseudo-effective iff  it is  effective. \qed
\end{cor}

\section{Maps of Calabi-Yau varieties}

Every variety has lots of different dominant rational maps
to projective spaces, but usually very few
dominant rational maps whose targets are not unirational.
The main result of this section proves a  
version of this for Calabi-Yau varieties.

\begin{thm} \label{iit.str.lem} Let $X$ be a projective Calabi-Yau variety
and $g:X\map Z$ a  dominant rational map
such that $Z$ is not uniruled.

Then there are
\begin{enumerate}
\item  a finite Calabi-Yau cover $h_X:\tilde X\to X$,
\item  an isomorphism $\tilde X\cong \tilde F\times \tilde Z$
where $\tilde F, \tilde Z$ are  Calabi-Yau varieties 
and $\pi_Z$ denotes the projection onto $\tilde Z$ , and
\item a generically finite map  $g_Z:\tilde Z\map Z$
\end{enumerate}
\noindent such that 
$g\circ h_X=g_Z\circ \pi_Z$.
\end{thm}

\begin{rem} 1. De-Qi Zhang pointed out to us that
in general $g_Z$ can not be chosen to be Galois
(contrary to our original claim). A simple example is given 
as follows. Let $A$ be an Abelian surface, $K=A/\{\pm 1\}$
the corresponding smooth Kummer surface and $m>1$. Then 
multiplication by $m$ on $A$ descends to a rational
map $K\map K$, but it is not Galois.

2. Standard methods of the Iitaka conjecture
(see especially \cite{kaw}) imply that for any dominant rational
map $X\map Z$,  either
 $\kappa(Z)= -\infty$
or  $\pi$  is an \'etale locally trivial fiber bundle
with Calabi--Yau fiber over an open subset of $Z$.
Furthermore,  \cite[Sec.8]{kaw} proves Theorem~\ref{iit.str.lem}
for the Albanese morphism.
The papers \cite{zhang1, zhang2} also contain related results and
techniques.

First  we explain how to modify the standard approach to the
Iitaka conjecture  to replace  $\kappa(Z)= -\infty$
with $Z$ uniruled. 

The remaining steps are more subtle since we have to construct
a suitable birational model of $Z$
and then  to extend the
fiber bundle structure from the open set to everywhere, at
least after a finite cover.
\end{rem}

Proof. As we noted before, there is a finite Calabi--Yau cover
$X'\to X$ such that $K_{X'}\sim 0$. 
In particular, $X'$ has canonical singularities.
We can further replace $Z$ by its normalization
in $\c(X)$. Thus we can assume to start with  that $K_X\sim 0$,
 $X$ has canonical singularities and $g$ has irreducible general fibers.

If $g$ is not a morphism along the closure of the general fiber of $g$
then $Z$ is uniruled. If $X$ is smooth, this is proved in
\cite[VI.1.9]{rc-book}; the general case follows from
\cite{ha-mc1}. Thus there  are open subsets $X^*\subset X$
and $Z^*\subset Z$ such that
$g:X^*\to Z^*$ is proper. We are free to shrink $Z^*$
in the sequel if necessary.

Let us look at a general fiber  $F\subset X$ of $g$.
It is a local complete intersection subvariety  whose
normal bundle is trivial. So, by the adjunction formula,
the  canonical class $K_F$ is
also trivial. $F$ has canonical singularities by
\cite[1.13]{reid}.

Choose smooth birational models 
$\sigma:X'\to X$ and $Z'\to Z$
 such that the corresponding $g':X'\to Z'$
is a morphism which   is smooth over the complement of
a simple normal crossing divisor  $B'\subset Z'$.
We can also assume that the image of every divisor
in $X\setminus X^*$ is a divisor in $Z'$.
Thus we can choose smooth open subvarieties 
$X^0\subset X$ and $Z^*\subset Z^0\subset Z'$ such that
\begin{enumerate}
\item $X\setminus X^0$ has codimension $\geq 2$ in $X$, and
\item $g_0:=g'|_{X^0}:X^0\to Z^0$ is flat and surjective (but not proper). 
\end{enumerate} 

The proof proceeds in 3 steps.

First, we show that 
 $\omega_{X'/Z'}|_{X^0}$ is the pull back of a line bundle
$L$ from $Z^0$ which is $\q$-linearly equivalent to 0.

Second we prove that there is an \'etale cover $Z^1\to Z^0$ 
and a birational map $Z^1\times F\map X^0\times_{Z^0}Z^1$
which is an isomorphism in codimension 1.

Third, we show that if $\tilde X\to X$ is the corresponding cover, then
  $\tilde X$ is a product of
two Calabi-Yau varieties, as expected.

In order to start with Step 1,
we need the following result about algebraic fiberspaces
for which  we could not find a convenient simple reference.

\begin{prop} \label{fib.space.prop}
 Notation as above. Then $g'_*\omega_{X'/Z'}$
is a line bundle and 
one can write the corresponding Cartier divisor
$$
\operatorname{divisor\ class\ of\ }
\bigl(g'_*\omega_{X'/Z'}\bigr)\simq J_g+B_g,
$$
where 
\begin{enumerate}
\item $B_g$ is an effective $\q$-divisor supported on $B'$, and
\item $J_g$ is a nef $\q$-divisor such that
\begin{enumerate}
\item either $J_g\simq 0$ and $g$ is an \'etale locally trivial fiber bundle
over some open set of $Z^*$,
\item or $(J_g\cdot C)>0$ for  every irreducible curve $C\subset Z'$
which is not contained in $B'$ and is not 
tangent to a certain foliation of $Z^*$.
\end{enumerate}
\end{enumerate}
\end{prop}

Proof. Over the open set where $g'$ is smooth, the results of
\cite[Thm.5.2]{griff} endow $g'_*\omega_{X'/Z'}$ with a Hermitian
metric whose curvature is semipositive.  This metric degenerates
along $B'$ but this degeneration is understood
\cite{fuj, kaw-abvar}, giving the decomposition
$J_g+B_g$ where $J_g$ is the curvature term and $B_g$ comes from
the singularities of the metric along $B'$.

Set $d=\dim F$, where $F$ is a general fiber of $g$.
If $F$ is smooth, 
we can assume that $F$ is also a fiber of $g'$. In this case
the curvature is flat in the directions
corresponding to the (left) kernel of the Kodaira--Spencer map
$$
H^1(F,T_F)\times H^0(F, \Omega_F^{d})\to H^1(F, \Omega_F^{d-1}).
$$
If $\Omega_F^{d}\cong \o_F$ then this is identified with the
Serre duality isomorphism
$$
  \bigl(H^{d-1}(F, \Omega_F^1)\bigr)^*\cong   H^1(F, \Omega_F^{d-1}),
$$
hence the above (left) kernel is zero. 
Thus $(J_g\cdot C)=0$ iff 
the deformation of the fibers $g^{-1}(C)\to C$ is trivial to first order
over every  point of $C\setminus B'$. This holds iff
the fibers of $g$ are all
isomorphic to each other over $C\setminus B'$.

The corresponding result for the case when $F$
has canonical singularities is worked out in \cite[Sec.6]{kaw}.\qed
\medskip

Let us now look at the natural map 
$$
g'^*\bigl(\omega_{Z'}\otimes g'_*\omega_{X'/Z'}\bigr)\to \omega_{X'}
$$
which is an isomorphism generically along $F$, thus nonzero.
Hence there is an effective divisor $D_1$ such that
$$
g'^*\bigl(\omega_{Z'}\otimes g'_*\omega_{X'/Z'}\bigr)\cong \omega_{X'}(-D_1).
$$
Write $\omega_{X'}\cong
\sigma^*\omega_X(D_2)\cong \o_{X'}(D_2)$,
where $D_2$ is $\sigma$-exceptional.

Let $C\subset X^0$ be a  general complete intersection 
curve.
 Then $\sigma^{-1}$ is
defined along $C$ and setting $C':=\sigma^{-1}(C)$ we get that
$$
\deg_{C'}g'^*\bigl(\omega_{Z'}\otimes g'_*\omega_{X'/Z'}\bigr)=
\deg_{C'}\omega_{X'}(-D_1)=\bigl(C'\cdot (D_2-D_1)\bigr)=
-\bigl(C'\cdot D_1\bigr).
$$
By the projection formula this implies that
$$
\bigl(g_0(C)\cdot K_{Z'}\bigr)+ \bigl(g_0(C)\cdot J_g\bigr)
+ \bigl(g_0(C)\cdot B_g\bigr)+ (C\cdot \sigma_*(D_1))=0.
\eqno{(*)}
$$
If $\bigl(g_0(C)\cdot K_{Z'}\bigr)<0$ then
$Z'$ is uniruled by the Miyaoka-Mori criterion \cite{mi-mo},
contrary to our assumptions. Thus all the terms on the left
hand side are nonnegative, hence they are all zero.

Since $C$ is a general curve, it can be chosen to be not tangent to
any given foliation. Therefore 
 $J_g$ is torsion in $\pic(Z')$
and $X^*\to Z^*$ is an \'etale locally trivial fiber bundle
for a suitable $Z^*$.
A general complete intersection curve 
 intersects every divisor in $X$, thus $(C\cdot \sigma_*(D_1))=0$
implies that  $\sigma_*(D_1)=0$, that is, $D_1$ is $\sigma$-exceptional.

Similarly,  $g_0(C)$ intersects every irreducible component
of $Z^0\setminus Z^*$. Thus $(g_0(C)\cdot B_g)=0$ implies that
$B_g|_{Z^0}=0$. These together imply that
 $L:=\bigl(g'_*\omega_{X'/Z'}\bigr)|_{Z^0}$
is $\q$-linearly equivalent to 0
and $\omega_{X'/Z'}|_{X^0}\cong g_0^*L$. This completes the first step.

Now to Step 2. 
Apply  Lemma~\ref{triv.after.et.cov} to $X^*\to Z^*$.
We get a  finite cover
$\pi:Z^1\to Z^0$ such that $X'\times_{Z^0}Z^1$ is birational to
$F\times Z^1$. By shrinking $Z^0$, we may assume that $Z^1$ is also smooth.
Eventually we prove that $\pi$ is \'etale over $Z^0$, but
for now we allow ramification over $Z^0\setminus Z^*$.

Let $n:X^1\to X'\times_{Z^0}Z^1$ be the normalization.
We  compare the relative dualizing sheaves
$$
\omega_{F\times Z^1/Z^1}\cong \o_{F\times Z^1}
\qtq{and} \omega_{X^1/Z^1}.
$$
Let $X^1\stackrel{u}{\leftarrow} Y\stackrel{v}{\to} F\times Z^1$ 
be a common resolution.
We can then write
$$
\omega_{Y/Z^1}\cong v^*\omega_{F\times Z^1/Z^1}(E_v)
\cong  \o_{Y}(E_v)
$$
for some divisor  $E_v$ supported on $\ex(v)$ and also
$$
\omega_{Y/Z^1}\cong u^* \omega_{X^1/Z^1}(E'_u)
\cong (g'\circ n\circ u)^*L(E_u)
$$
for some divisors $E'_u,E_u$.
Since $g^*L|_{X^0}$ is $\q$-linearly equivalent to zero,
we conclude from these that 
$$
u_*(E_u-E_v) |_{X^0}\simq 0.
$$
Next we get some information about $E_u$ and $E_v$.
Since $F\times Z^1$ has canonical singularities, 
every irreducible component of  $E_v$  is effective.
Furthermore, an irreducible component of  $\ex(v)$ 
appears with positive coefficient in $E_v$
unless it dominates $Z^0$. Thus we see that 
$u_*(E_v) |_{X^0}$ is supported in $X^0\setminus X^*$
and an  irreducible component of $X^0\setminus X^*$
appears with positive coefficient in $u_*(E_v)$, unless
$v\circ u^{-1}$ is a local isomorphism over  its generic point.

On the other hand, since $\omega_{X'/Z'}|_{X^0}$ is the pull back of $L$,
  $\omega_{X'\times_{Z^0}Z^1/Z^1}|_{X^0}$ is 
the pull back of  $\pi^*L$. 
As we normalize, we subtract divisors corresponding to the
nonnormal locus. 
Every other irreducible component of $E_u$ is $u$-exceptional,
hence gets killed by $u_*$.
Thus we obtain that 
$u_*(E_u) |_{X^0}$ is also
contained  in $X^0\setminus X^*$ 
and an  irreducible component of $X^0\setminus X^*$
either appears with negative coefficient or 
 $X'\times_{Z'}Z^1$ is normal over that component
and the coefficient is 0. 

Thus $u_*(E_u-E_v) |_{X^0}$ is a nonpositive linear combination of
the irreducible components of $X^0\setminus X^*$
and it is also $\q$-lineraly equivalent to 0.
Since $X\setminus X^0$ has codimension $\geq 2$,
we conclude that $u_*(E_u-E_v) |_{X^0}=0$.
That is, $X^0\times_{Z^0}Z^1$ is normal 
in codimension 1
and isomorphic to $F\times Z^1$, again only in codimension 1.

We may as well assume that $Z^1\to Z^0$ is Galois with group $G$.
We then have a corresponding $G$-action on
$F\times Z^1$ for a suitable $G$-action on $F$.
By taking the quotient, we obtain a birational 
map 
$$
\phi:W:=(F\times Z^1)/G\to Z^0 \map X^0
$$
which is an isomorphism in codimension 1. In particular,
$$
\omega_{W/Z^0}\cong \phi^*\omega_{X^0/Z^0}\cong g_W^*L
$$
where $g_W:W\to Z^0$ is the quotient of the projection map to $Z^1$.
Using  (\ref{quots.of.triv}) we conclude that
 $g_W:W\to Z^0$ is in fact an  \'etale locally trivial fiber bundle
with fiber $F$, at least outside a codimension 2 set.

Then  Lemma~\ref{triv.after.et.cov}
shows that $Z^1\to Z^0$ is also  \'etale at every
generic point of $Z^0\setminus Z^*$, thus it is
a finite \'etale cover.

Let now $\tilde X\to X$ be the normalization
of $X$ in the function field of $X^1$. 
Since $Z^1\to Z^0$ is \'etale, we see that $\tilde X\to X$
is \'etale over $X^0$. Thus $\tilde X\to X$
is  \'etale  outside a codimension $\geq 2$ set.
In particular, $\tilde X$ is a Calabi-Yau variety.

Furthermore, 
the  birational map $\phi:\tilde X\to F\times Z^1$
 is an open embedding outside a codimension $\geq 2$ set.
That is, $\phi$ does not contract any divisor. This completes Step 2.

By Proposition~\ref{bir.to.prod}, $\tilde X$ is itself a product 
$\tilde F\times \tilde Z$.
Note that  $F=\tilde F$ since $\phi$ is an isomorphism
along $F\cong \pi^{-1}(z)$ for $z\in Z^*$.
Thus  $\tilde X\cong F\times \tilde Z$
and there is a generically finite map $\tilde Z\map Z$.
\qed

\begin{lem}\label{triv.after.et.cov}
Let $f:U\to V$ be a projective morphism between normal varieties. $V$ smooth.
Assume that $f$ is an \'etale locally trivial fiber bundle with typical
fiber $F$ which is a Calabi-Yau variety. Then there is a finite \'etale
cover $V'\to V$ such that the pull back $U\times_VV'\to V'$ is
globally trivial. Moreover, we can choose $V'\to V$ such that
its generic fiber depends only on the generic fiber of $f$.
\end{lem}

Proof. Let $H$ be an ample divisor on $U$.
Let $\pi:\isom(F\times V, U, H)\to V$ denote the $V$-scheme 
parametrizing $V$-isomorphisms $\phi:F\times V\to  U$ such that
 $\phi^*H$ is numerically equivalent  to $H_F$.
The  fiber of $\isom(F\times V, U, H)\to V$
over $v\in V$ is the set of isomorphisms $\phi:F\to U_v$
such that $\phi^*(H|_{U_v})$ is numerically equivalent  to $H_F$.

Note that $\isom(F\times V, U, H)\to V$
is an  \'etale locally trivial fiber bundle with typical
fiber $\aut_H(F)$. 
Any \'etale multisection of $\pi$ gives a required
\'etale cover $V'\to V$.

Thus we need to find an \'etale multisection of
a projective group scheme (in characteristic 0).
The Stein factorization of $\pi$ 
gives an \'etale cover $V_1\to V$ and if we pull back everything to
$V_1$, then there is a 
well defined identity component.
Thus we are reduced to the case when $\pi:I\to V$
 is torsor over an Abelian scheme  $A\to V$.

Let $I_g$ be the generic fiber and let $P\in I_g$ be a point of degree $d$.
Let $S_d\subset I_g$ be the set of geometric points
$p$ such that $dp-P= 0\in A_g$. Then $S_d$ is defined over 
$k(V)$
and it is a principal homogeneous space over the subgroup of
 $d$-torsion points of $A_g$. We claim that the closure of $S_d$ in $I$
 is finite and 
\'etale over $V$. Indeed, it is finite over
codimension 1 points and also \'etale over
codimension 1 points since the limit of  nonzero $d$-torsion points 
cannot be zero. Thus it is also \'etale over all points by
the purity of branch loci.\qed
\medskip

The $K_X=0$ of the following lemma is proved in
\cite[Thm.2]{pet-min}. 

\begin{prop}\label{bir.to.prod}
 Let $X, U,V$ be normal projective varieties. Assume that
$X$ has rational singularities.
Let $\phi:X\map U\times V$ be a birational map which does not contract
any divisor. Then there are  normal projective varieties
$U'$ birational to $U$ and $V'$ birational to $V$ such that
$X\cong U'\times V'$.
\end{prop}

Proof. We can replace $U,V$ by resolutions, thus we may assume that 
they are smooth.

Let $X\stackrel{p}{\leftarrow} Y \stackrel{q}{\to} U\times V$
be a factorization of $g$. By assumption, $\ex(q)\subset \ex(p)$.
Let $H$ be a very ample divisor on $X$ and $\phi_*H=q_*p^*H$
its birational transform.
Then
$| q_*p^*H|=| p^*H+m\ex(q)|$ for $m\gg 1$.
On the other hand,
$$
|H|=| p^*H|\subset | p^*H+m\ex(q)|\subset  | p^*H+m\ex(p)|=|H|.
$$ 
Thus $|H|=|\phi_*H|$.

Assume that
there are divisors $H_U$ on $U$ and $H_V$ on $V$ such that
$\phi_*H\sim \pi_U^*H_U+\pi_V^*H_V$. Then 
$$
H^0(U\times V, \o_{U\times V}(\phi_*H))=
H^0(U,\o_U(H_U))\otimes H^0(V,\o_V(H_V)).
$$
Since $X$ is the closure of the image of $U\times V$
under the linear system $|\phi_*H|$, we see that
$X\cong U'\times V'$ where  $U'$ is the image of $U$
under the linear system $|H_U|$ and $V'$ is the image of $V$
under the linear system $|H_V|$.

 If $H^1(U,\o_U)=0$, then
$\pic(U\times V)=\pi_U^*\pic(U)+ \pi_V^*\pic(V)$, and the
we are done.
In general, however, 
$\pic(U\times V)\supsetneq \pi_U^*\pic(U)+ \pi_V^*\pic(V)$,
and we have to change $H$.

 Fix
points $u\in U$ and $v\in V$ and let
$D_U:=\phi_*H|_{U\times\{v\}}$ and $D_V:=\phi_*H|_{\{u\}\times V}$.
Set $D':=\phi_*H-\pi_U^*D_U-\pi_V^*D_V$.
Then $D'$ restricted to any $U\times\{v'\}$ is in $\pic^0(U)$ and
$D'$ restricted to any $\{u'\}\times V$ is in $\pic^0(V)$.
Thus there is a divisor $B$ on $\Alb(U\times V)$
such that  $D'=\alb_{U\times V}^*B$, where, for a variety $Z$,
$\alb_Z:Z\to \Alb(Z)$ denotes the Albanese map.

Choose divisors $B_U$ on $\Alb(U)$ and  $B_V$ on $\Alb(V)$
such that $\pi_U^*B_U+\pi_V^*B_V-B$ is very ample, where,
somewhat sloppily, $\pi_U,\pi_V$ also denote the coordinate projections of
$\Alb(U\times V)$.

Since $X$ has rational singularities,
$\Alb(X)=\Alb(U\times V)$.
Replace $H$ by
$$
H^*:=H+\alb_X^*(\pi_U^*B_U+\pi_V^*B_V-B).
$$
Then
$$
\begin{array}{rcl}
\phi_*H^*&=&\phi_*H+\alb_{U\times V}^*(\pi_U^*B_U+\pi_V^*B_V-B)\\
&=& \pi_U^*(H_U+\alb_U^*B_U)+\pi_V^*(H_V+\alb_V^*B_V).
\end{array}
$$
Since $H^*$ is again very ample, we are done.
\qed

\medskip

\begin{say}[Quotients of trivial families]\label{quots.of.triv}
We consider families $X\to C$ over a smooth pointed curve germ
$0\in C$ such that after a finite base change and normalization
we get a trivial family. This means that we start with a trivial family
$F\times D$ over a disc $D$, an automorphism $\tau$ of $F$ of order dividing
$m$
and take the quotient
$X:=(F\times D)/(\tau, e(1/m))$.

If the order of $\tau$ is less than $m$ then there is a subgroup
which acts trivially on $F$ and the quotient is again a trivial family.
Thus we may assume  that the order of $\tau$ is precisely $m$.

Fix a top form $\omega$ on $F$.
Pulling back by $\tau$ gives an isomorphism
$\omega=\eta  \tau^*\omega$ for some 
$m$th root of unity $\eta$. If $\eta\neq 1$ then 
on the quotient family the monordomy around $0\in C$ has finite order
 $\neq 1$, and the boundary term $B$ in Proposition~\ref{fib.space.prop} is nonzero.

Finally, if $\omega= \tau^*\omega$ 
then $F_0:=F/(\tau)$ also has trivial canonical class. Thus by the
adjunction formula we see that
$\omega_{X/C}\cong 
\o_{X}((m-1)F_0)$ 
 is not trivial.
\end{say}

Next we apply Theorem~\ref{iit.str.lem}
to study those quotients of Calabi-Yau varieties
which are  uniruled but not
rationally connected.
Let us see first some examples of how this can happen.
Then we see that these trivial examples exhaust all
possibilities.

\begin{exmp} Let $\Pi:X'\to Z'$ be an \'etale locally trivial
fiber bundle whose base $Z'$ and typical fiber $F'$ are both 
 projective Calabi-Yau varieties. Then $X'$ is also a  projective
Calabi-Yau variety. Let $G'$ be a finite group acting on $X'$
and assume that $\Pi$ is $G'$-equivariant. Assume that
\begin{enumerate}
\item $\kappa(Z'/G')=0$, and
\item for general $z\in Z'$, the quotient
$\Pi^{-1}(z)/\stab_zG'$ is rationally connected.
\end {enumerate}
Then $\Pi/G':X'/G'\to Z'/G'$ is the MRC-fibration of
$X'/G'$.

More generally, let $H\subset G'$ be a normal subgroup
such that $X:=X'/H$ is a Calabi-Yau variety
and set $G:=G'/H$.
Then $X/G\cong X'/G'$ and its MRC-fibration
is given by $\Pi/G':X'/G'\to Z'/G'$.
\end{exmp}

\begin{thm} \label{cy.nonsimple.thm} 
Let $X$ be a  projective Calabi--Yau variety
and $G$ a finite group acting on $X$. 
Assume that $X/G$ is  uniruled but not rationally connected.
Let $\pi:X/G\map Z$ be the MRC fibration.
Then there is
\begin{enumerate}
\item a finite, Calabi--Yau,  Galois cover $X'\to X$,
\item a proper morphism $\Pi: X'\to  Z'$ which is
 an \'etale locally trivial
fiber bundle whose base $ Z'$ and typical fiber $ F'$ are both 
 projective Calabi-Yau varieties, and
\item a group $ G'$ acting on $ X'$
where  $ G\supset \gal(X'/X)$ and
$ G'/\gal(X'/X)=G$,
\end{enumerate}
such that
  $\Pi/G':X'/G'\to  Z'/G'$
 is birational to the MRC fibration
$\pi:X/G\map Z$.
\end{thm}

Proof.   
Let $X/G\map Z$ be the MRC fibration and let
$\pi:X\to X/G \map Z$ be the composite.
$Z$ is not uniruled by
\cite{ghs} thus Theorem~\ref{iit.str.lem} applies and we get
a direct product $F\times Z$ mapping to
$X$. Since both $X$ and $F\times Z$ have trivial canonical class,
$F\times Z\to X$ is \'etale in codimension 1.

 In order to lift the $G$-action from $X$ to a cover, we need
to take the Galois closure of  $F\times Z\to X/G$.
Let $G'$ be its Galois group.
This replaces $F\times Z$ with a finite cover which is 
\'etale in codimension 1. The latter
need not be globally a product, only  \'etale locally a product.\qed

\begin{cor} \label{abvar.max.quot}
Let $A$ be an Abelian variety and  $G$ a
 finite group acting on $A$.
There is a unique maximal $G$-equivariant quotient $A\to B$
such that
$A/G\to B/G$  is the MRC quotient.
\end{cor}

Proof. Let $\Pi:A^0\to Z^0$ be the quotient
constructed in Theorem~\ref{cy.nonsimple.thm}. Its fibers $F_z$ are smooth 
subvarieties of $A$ with trivial canonical class.
Thus each  $F_z$ is a translation of a fixed
 Abelian subvariety $C\subset A$ (cf.\ \cite[4.14]{gr-ha}).
Set $B=A/C$.\qed

\begin{defn} Let $G$ be a finite group acting
on a vector space $V$. Let $G^{RT}<G$ be the subgroup generated by
all elements of age $<1$ and $V^{RT}$ the complement of the
fixed space of $G^{RT}$.
\end{defn}

\begin{defn} Let $A$ be an Abelian variety and $G$ a finite group acting
on $A$. For every $x\in A$, let $G_x:=\stab(x)< G$ denote the stabilizer
and $i_x:A\to A$ the translation by $x$.
Consider the action of $G_x$ on $T_xA$, the tangent space of $A$ at $x$.
Let 
$G_x^{RT}$ and $(T_xA)^{RT}$ be as above.
Note that $(T_xA)^{RT}$ is the tangent space of a translate of
an Abelian subvariety $A_x\subset A$
since it is the intersection of the kernels of the endomorphisms
$g-{\mathbf 1}_A$ for $g\in G_x^{RT}$.
Set
$$
G^{RT}:=\bigl\langle G_x^{RT}: x\in A\bigr\rangle\qtq{and}
(TA)^{RT}:=\bigl\langle i_x^*(T_xA)^{RT}: x\in A\bigr\rangle.
$$
Then $(TA)^{RT}$ is  the tangent space of the Abelian subvariety
generated by the $A_x$. Denote it by $A^{RT}_1$.

The group $G/G^{RT}$ acts on the quotient Abelian variety 
$q_1:A\to A/A^{RT}_1$.  If $q_i:A\to A/A^{RT}_i$
is already defined, set 
$$
A^{RT}_{i+1}:=q_i^{-1}\left((A/A^{RT}_i)^{RT}_1\right)
$$
and let $q_{i+1}:A\to A/A^{RT}_{i+1}$ be the quotient map.
The increasing sequence of Abelian subvarieties 
$A^{RT}_1\subset A^{RT}_2 \subset \cdots$ eventually
stabilizes to $A^{RT}_{stab}\subset A$.
\end{defn}

\begin{cor} 
\label{filtered.rc}
Let $A$ be an Abelian variety and $G$ a finite group acting
on $A$. Then
\begin{enumerate}
\item $\kappa(A/G)=0$ iff $G^{RT}=\{1\}$, and
\item $A/G$ is rationally connected  iff $A^{RT}_{stab}=A$.
\end{enumerate}
\end{cor}

\section{Basic non-Reid-Tai pairs}

Our goal in this section is to classify basic non-RT pairs.

There is a basic dichotomy:

\begin{prop}
\label{Prod}
If $(G,V)$ is a basic non-RT pair, then either $(G,V)$ is
primitive or $G$ respects a decomposition of
$V$ as a direct sum of lines: 
$$V = L_1\oplus\cdots\oplus L_n.$$
In the latter case, the 
homomorphism $\phi\colon G\to S_n$ given by the permutation action of $G$ on $\{L_1,\ldots,L_n\}$ is
surjective, and every exceptional 
element in $G$ maps to a transposition.
\end{prop}

\begin{proof}
Suppose that $G$ respects the decomposition $V\cong V_1\oplus\cdots\oplus V_m$ for 
some $m\ge 2$.   If there is more than one such decomposition, we choose one so that
$m\ge 2$ is minimal.  By irreducibility, $G$ acts transitively on the set of $V_i$.  
As the conjugacy class of
some exceptional element $g$ generates $G$, it follows that $g$ permutes the $V_i$
non-trivially.  Suppose that for $2\le k\le m$, we have
$$g(V_1) = V_2,\,g(V_2)=V_3,\, g(V_k)=V_1.$$
Then $g$ and $e(1/k) g$ are isospectral on $V_1\oplus\cdots\oplus V_k$.
Thus the eigenvalues of $g$ constitute a union of $\dim V_1$ cosets of the cyclic group
$\langle e(1/k)\rangle$.  Such a union of cosets can satisfy the Reid-Tai condition
only if $k=2$ and $\dim V_1 = 1$, and then $g$ must stabilize $V_i$ for every $i\ge 3$.
Thus $g$ induces a transposition on the $V_i$, each of which must be of dimension $1$.
A transitive subgroup of $S_m$ which contains a transposition must be of the form 
$S_a^b\rtimes T$, where $ab=m$, $T$ is a transitive subgroup of $S_b$, and $a\ge 2$.  
It corresponds to a decomposition of the set of factors $V_i$ into $b$ sets of cardinality $a$.
If $W_j$ denote the direct sums of the $V_i$ within each of the $a$-element sets of this
partition, it follows that $G$ respects this coarser decomposition, contrary to the assumption
that $m$ is minimal.
\end{proof}

To analyze the primitive case, it is useful to quantify how far a unitary operator is from the identity.

\begin{defn}
Let $H$ be a Hilbert space, $T$ a unitary operator on $H$, and $B$ an orthonormal basis of
$H$.  The \emph{deviation} of $T$ with respect to $B$ is given by
$$d(T,B):=\sum_{b\in B} \Vert T(b)-b\Vert.$$
The \emph{deviation} of $T$ is 
$$d(T) := \inf_B d(T,B),$$
as $B$ ranges over all orthonormal bases.  If $d(T) < \infty$, we say $T$ has \emph{finite
deviation}. 
\end{defn}

As the arc of a circle cut off by a chord is always longer than the chord, if $H$ is a finite-dimensional Hilbert space and $g\colon H\to H$ a unitary operator of finite order
satisfying the Reid-Tai condition, we have $d(g) < 2\pi$.  This is the primary motivation
for our definition of deviation.

For any space $H$, unitary operator $T$, basis $B$, and real number $x>0$, we define
$$S(T,B,x) = \{b\in B\mid \Vert T(b)-b\Vert\ge x\}.$$
If $1_I$ denotes the characteristic function of the interval $I$, 
\begin{align*}\int_0^\infty |S(T,B,x)| dx &= \int_0^\infty \sum_{b\in B} 1_{[0,\Vert T(b)-b\Vert]}dx \\
&= \sum_{b\in B} \int_0^\infty 1_{[0,\Vert T(b)-b\Vert]}dx = \sum_{b\in B}\Vert T(b)-b\Vert
= d(T,B).\\
\end{align*}

It is obvious that deviation is symmetric in the sense that
$d(T)=d(T^{-1})$.  Next we prove a lemma relating $d$ to multiplication in
the unitary group.

\begin{prop}
\label{product.bound}
If $T_1,T_2,\ldots T_n$ are unitary operators of finite deviation on a Hilbert space $H$, then
$$d(T_1 T_2\cdots T_n)\le n(d(T_1)+d(T_2)+\cdots d(T_n)).$$
\end{prop}

\begin{proof}
Let $B_1, B_2,\ldots, B_n$ denote orthonormal bases of $H$.  We claim that there exists 
an orthonormal basis $B$ such that for all $x>0$,
\begin{equation}
\label{prod-est}
|S(T_1T_2\cdots T_n,B,nx)|\le |S(T_1,B_1,x)|+|S(T_2,B_2,x)|+\cdots |S(T_n,B_n,x)|.
\end{equation}
This claim implies the proposition, by integrating over $x$.

Given $T_i, B_i$ and $x>0$, we define
$$V_x = \Span \bigcup_{i=1}^n S(T_i,B_i,x).$$
As all $T_i$ are of finite deviation, the set of ``jumps'' ($x$ such that
$V_x$ is not constant in a neighborhood of $x$) is discrete in $(0,\infty)$.
Arranging them in reverse order, we see that there exists a (possibly infinite) increasing
chain of finite-dimensional subspaces $W_i$ of $H$ such that each $V_x$ is equal to one
of the $W_i$.  We choose $B$ to be any orthonormal basis adapted to 
$W_1\subset W_2\subset \cdots$ in the sense that $B\cap W_i$ is an orthonormal basis
of $W_i$ for all $i$.

For all $b\in B$, by the triangle inequality,
$$\Vert T_1 T_2\cdots T_n(b) - b\Vert  \le 
\sum_{i=1}^n \Vert T_1 T_2\cdots T_{i-1}(T_i(b) - b)\Vert
= \sum_{i=1}^n \Vert T_i(b)-b\Vert.$$
If $b\notin V_x$, then $b$ is orthogonal to every element of $S(T_i,B_i,x)$ for
$i=1,\ldots,n$, and therefore, $\Vert T_i(b) - b\Vert \le x$ for all $i$.  It follows that
$$\Vert T_1 T_2\cdots T_n(b) - b\Vert  \le nx,$$
or $b\notin S(T_1T_2\cdots T_n,B,nx)$.  This implies (\ref{prod-est}).

\end{proof}

\begin{prop}
\label{Comm}
If $T_1$ and $T_2$ are operators on a Hilbert space $H$ such that $T_1$ is of bounded
deviation, then 
$$d(T_1^{-1} T_2^{-1} T_1 T_2) \le 4 d(T_1).$$
\end{prop}

\begin{proof}
As $d(T_1^{-1}) = d(T_1) = d(T_2^{-1} T_1 T_2)$, the proposition follows from Proposition~\ref{product.bound}.
\end{proof}

\begin{lemma}
\label{big.deviation}
Let $G$ be a compact group and $(\rho,V)$ a non-trivial representation of $G$ such that $V^G = (0)$.
Then there exists $g\in G$ with $d(g)\ge \dim V$.
\end{lemma}

\begin{proof}
As $V$ has no $G$-invariants,
$$\int_G \tr(\rho(g))\,dg = 0,$$
there exists $g\in G$ with $\Re(\rho(g))\le 0$.  If $d(g) < \dim V$, there exists 
an orthonormal basis $B$ of $V$ such that $\sum_{b\in B} |g(b)-b| < \dim V$.
If $a_{ij}$ is the matrix of $\rho(g)$ with respect to such a basis, 
$\Re(\sum_i a_{ii}) < 0$, so 
$$\sum_{b\in B} \Vert b-g(b)\Vert > \sum_i |1-a_{ii}| \ge \sum_i (1-\Re(a_{ii})) > \dim V,$$
which gives a contradiction.

\end{proof}

\begin{lemma}
\label{bounded.cn}
Let $(G,V)$ be a finite group and a representation such that $V^G = (0)$.  
Suppose that for some integer $k$,
every element of $G$ can be written as a product of at most $k$ conjugates of $g$.
Then $d(g) \ge \frac{\dim V}{k^2}$.
\end{lemma}

\begin{proof}
This is an immediate consequence of Proposition~\ref{product.bound}
and Lemma~\ref{big.deviation}.
\end{proof}

We recall that a \emph{characteristically simple} group $G$ is isomorphic to a group of the
form $K^r$, where $K$ is a (possibly abelian) finite simple group.

\begin{prop}
\label{log.cn}
There exists a constant $C$ such that if $H$ is a perfect central extension of a characteristically
simple group $K^r$ and the conjugacy class of $h\in H$ generates the whole group, then
every element in $H$ is the product of no more than $C\log|H|$ conjugates of $h$.
\end{prop}

\begin{proof}
If $H$ is perfect, then $K^r$ is perfect, so $K$ is a non-abelian finite simple group and therefore
perfect.  Let $\tilde K$ denote the universal central extension of $K$, so $\tilde K^r$ is
the universal central extension of $K^r$ as well as of $H$.   As the only subgroup of $\tilde K^r$
mapping onto $H$ is $\tilde K^r$, it  suffices to prove that
if $\tilde h = (x_1,\ldots,x_r)\in \tilde K^r$ has the property that its 
conjugacy class generates $\tilde K^r$,
then every element of $\tilde K^r$ is the product of at most $C\log |H|$ conjugates of $\tilde h$.
For any $t$ from $1$ to $r$, we can choose $(1,\ldots,1,y,1,\ldots,1)\in \tilde K^r$
(with $y$ in the $t$th coordinate) whose commutator with $\tilde h$
is an element $(1,\ldots,1,z,1\ldots,1)$ not in the center of $\tilde K^r$.
This element is the product of two conjugates of $\tilde h$.  If we can find an absolute constant
$A$ such that for every finite simple group $K$ and every non-central element $z\in \tilde K$,
every element of $\tilde K$ is the product of at most $A\log |K|$ conjugates of $z$, 
the proposition holds with $C = 2A$.

To prove the existence of $A$, we note that we may assume that $K$ has order greater than any specified constant.  In particular, we may assume that $K$ is either an alternating group 
$A_m$, $m\ge 8$, or a group of Lie type, and in the latter case,
if $K\neq \tilde K$, we can identify $\tilde K$ with
the group of points of a simply connected, almost simple algebraic group over a finite field.
It is known that the \emph{covering number} of $A_m$ is
$\lfloor m/2\rfloor$ (see, e.g., \cite{ArHe}), so every element of the group can be written as a product of 
$\le m/2$ elements belonging to any given non-trivial conjugacy class $X$.  The universal central extension of $A_m$ is of order $m!$, and at least half of those elements are products of
$\le m/2$ elements in any fixed conjugacy class $X$, so all of the elements are products of 
$\le m < \log m!/2$ elements of $X$.  For the groups of Lie type, we have an upper bound linear in
the absolute rank (\cite{EGH},\cite{LaLi}) and therefore sublogarithmic in order.

\end{proof}

\begin{lemma}
\label{Jordan}
For every integer $n>0$, there are only finitely many classes of primitive 
finite subgroups $G\subset \GL_n(\C)$ up to projective equivalence.
\end{lemma}

\begin{proof}
As $G$ is primitive, an normal abelian subgroup of $G$ lies in the center of $\GL_n(\C)$.
By Jordan's theorem, $G$ has a normal abelian subgroup whose index can be bounded in terms 
of $n$.  Thus the image of $G$ in $\PGL_n(\C)$ is bounded in terms of $n$.  For each isomorphism
class of finite groups, there are only finitely many projective $n$-dimensional representations.
\end{proof}

\begin{say}[Proof of Theorem~\ref{nonRT.almost.class}]

First we assume $G$ stabilizes a set $\{L_1,\ldots,L_n\}$ of lines which give a
direct sum decomposition of $V$.  We have already seen that the resulting homomorphism
$\phi\colon G\to S_n$ is surjective.
Let $t\in G$ lie in $\ker\phi$, so
$t(v_i) = \lambda_i v_i$ for all $i$ and all $v_i\in L_i$.  The commutator of $t$ with any preimage of
the transposition $(i\, j)\in S_n$ gives an element of $\ker \phi$ that has eigenvalues
$\lambda_i/\lambda_j$, $\lambda_j/\lambda_i$ and $1$ (of multiplicity $n-2$).  The $G$-conjugacy
class of this element consists of all diagonal matrices with this multiset of eigenvalues.
Thus $\ker\phi$ contains $\ker \det_C\colon C^n\to C$, where $C$ is the group generated by all
ratios of eigenvalues of all elements of $\ker \phi$.  It follows that $\ker\phi$ is the product of
$\ker \det_C$ and a group of scalar matrices.  If we pass to $\PGL(V)$, therefore,
the image of $G$ is an extension of $C^{n-1}$ by $S_n$.

We claim that this extension is split if $n$ is sufficiently large.  To prove this, it suffices to prove
$H^2(S_n,C^{n-1}) = 0$.  This follows if we can show that $H^2(S_n,\Z^{n-1})=H^3(S_n,\Z^{n-1})=0$,
or equivalently, the sum-of-coordinate maps
$$H^i(S_n,\Z^n)\to H^i(S_n,\Z)$$ 
are isomorphisms for $i=1,2,3$, where $S_n$ acts on $\Z^n$ by permutations.  By Shapiro's lemma,
the composition of restriction and sum-of-coordinates gives an isomorphism
$H^i(S_n,\Z^n)\tilde{\to} H^i(S_{n-1},\Z)$, so we need to know that the restriction homomorphisms
$H^i(S_n,\Z)\to H^i(S_{n-1},\Z)$ are isomorphisms when $n$ is large compared to $i$, which
follows from \cite{nakaoka}.

Thus the image of $G$ in $\PGL(V)$ is $C^{n-1}\rtimes S_n$, which is the same as the image in 
$\PGL(V)$ of the imprimitive unitary reflection group $G(|C|,k,n)$, where $k$ is any divisor of $G$.

It remains to consider the primitive case.  Let $Z$ denote the center of 
$G$.  As $Z$ is abelian and has a faithful isotypic representation, it must be cyclic.
If $G$ is abelian, then $G=Z$, and we are done.  (This can be regarded as a subcase
of the case that $G$ stabilizes a decomposition of $V$ into lines.)  Otherwise,
let $\bar H\cong K^r$ denote a characteristically simple normal subgroup of $\bar G:=G/Z$,
where $K$ is a (possibly abelian) finite simple group and $r\ge 1$.  If $K$ is abelian, we 
let $H\subset G$ denote the inverse image of $\bar H$ in $G$.  
In the non-abelian case, we let $H$ denote the derived group of the inverse image
of $\bar H$ in $G$, which is perfect and again maps onto $\bar H$.
We know that $H$ is not contained in the center 
of $G$, so some inner automorphism of $G$ acts non-trivially on $H$.  It follows that conjugation
by $g$ acts non-trivially on $H$.  By Proposition~\ref{Comm}, there exists a non-trivial
element $h\in H$ with $d(h)< 8\pi$.

We consider five cases:
\begin{enumerate}
\item $H$ is abelian.
\item $\bar H$ is abelian but $H$ is not.
\item $K$ is a group of Lie type.
\item $K$ is an alternating group $A_m$, where $m$ is greater than a sufficiently large constant.
\item $K$ is non-abelian but not of type (3) or (4).
\end{enumerate}

We prove that case (4) leads to reflection groups, and all of the other cases contribute only
finitely many solutions.

If $H$ is abelian, then the restriction of $V$ to $H$ is isotypical and $H$ is central, contrary 
to the definition of $H$.

If $\bar H$ is abelian and $H$ is not, then $H$ is a central extension of a vector group and is therefore the product of its center $Z$ by an extraspecial $p$-group $H_p$
for some prime $p$.  The kernel $Z[p]$ of multiplication by $p$ on $Z$ is isomorphic to
$\Z/p\Z$, and the commutator map $\bar H\times \bar H \to Z[p]$ gives a non-degenerate pairing.  Therefore, conjugation by any element in $G\setminus H$ induces a non-trivial map on $\bar H$,
Thus we can take $h$ with $d(h) < 8\pi$ to lie outside the center of $H$.

The image of every element of $H$ in $\Aut(V)$ is the product of a scalar matrix and the
image of an element of $H_p$ in a direct sum of $m\ge 1$ copies of one of its faithful
irreducible representations.  By the Stone-von Neumann theorem, a faithful representation
of an extraspecial $p$-group is determined by a central character; its dimension is $p^n$,
where $|H_p| = p^{2n+1}$, and every non-central element has eigenvalues
$\omega,\omega e(1/p), \omega e(2/p),\ldots,\omega e(-1/p)$, each occurring with multiplicity $p^{n-1}$, where $\omega^{p^2}=1$.   As $h$ is a scalar multiple of an element with these eigenvalues,
$$8\pi > d(h) \ge 2\pi m p^{n-1}\frac{p(p-1)}{2p}\ge \frac{2\pi m p^n}4,$$
so $\dim V = m p^n < 16$.   By Lemma~\ref{Jordan}, there are only finitely many possibilities 
for $(G,V)$ up to projective equivalence.

In cases (3)--(5), $H$ is perfect.  If the conjugacy class of $h$ in $G$ does not generate $H$,
it generates a proper normal subgroup of $H$, which is a central extension of a
subgroup of $K^r$ which is again normal in $\bar G$.  
Such a subgroup is of the form $K^s$ for $s < r$.  Replacing $H$
if necessary by a smaller group, we may assume that the $G$-conjugacy class of $h$ generates
$H$.  By Proposition~\ref{log.cn}, every element of $H$ is the product of at most
$C\log|H|$ conjugates of $h$, and by Lemma~\ref{bounded.cn}, this implies 
\begin{equation}
\label{dim.lower.bound}
\dim V < 8\pi C^2 \log^2|H|.
\end{equation}

For case (3), we note that by \cite[Table 1]{sz}, a faithful irreducible projective
representation of a finite simple group $K$ which is not an alternating group
always has dimension at least $e^{c_1\sqrt{\log |K|}}$ for some positive absolute constant $c_1$.  A faithful irreducible representation of $K^r$ is the tensor power of
$r$ faithful irreducible representations of $K$, so its dimension is at least
$e^{c_1\sqrt{\log |\bar H|}} > e^{c_1\sqrt{\log |H|}/2}$.   For $|H|\gg 0$, this is in contradiction with (\ref{dim.lower.bound}).
Thus there are only finitely many possibilities for $H$ up to isomorphism, and this gives an upper bound for $\dim V$.  By Lemma~\ref{Jordan}, there are only finitely many possibilities for $(G,V)$
up to projective equivalence.

For case (4), we need to consider both ordinary representations of $A_m$
and spin representations (i.e., projective representations which do not lift to linear
representations).
By \cite{kt}, the minimal degree of a spin representation of $A_m$ grows faster than any
polynomial, in particular, faster than $m^{5/2}$.
To every irreducible linear representation of $A_m$ one can associate a partition $\lambda$ of $m$
for which the first part $\lambda_1$ is greater than or equal to the number of parts.
There may be one or two representations associated to $\lambda$, and in the latter case their 
degrees are equal and their direct sum is irreducible as an $S_m$-representation.
By \cite[2.1,2.4]{ls-fuchsian}, if $\lambda_1\le m-3$,
the degree of any $S_m$-representation associated to $\lambda$ is at least
$\binom {m-3}{3}$, so the degree of any $A_m$-representation 
associated to $\lambda$ is at least half of that.
Thus, for $m\gg 0$, the only faithful representations of $A_m$ which have degree less than 
$m^{5/2}$ are $V_{m-1,1}$, $V_{m-2,2}$, and $V_{m-2,1,1}$ of degrees
$m-1$, $\frac {(m-1)(m-2)}{2}$ and $\frac{m(m-3)}{2}$ respectively.
Now, $\log |H| \le \log (m!)^r < rm^{1.1}$ for $m\gg 0$, and the minimal degree
of any faithful representation of $H$ is at least $(m-1)^r \gg r^2 m^{2.2}$ for $r\ge 3$.
By (\ref{dim.lower.bound}), there are only four possibilities which need be considered.  
If $r=2$, then $H=\bar H = A_m^2$, and $V$ must be the
tensor product of two copies of $V_{m-1,1}$.  Otherwise $r=1$, $H=\bar H=A_m$, and $V$ is
$V_{m-1,1}$, $V_{m-2,2}$, or $V_{m-2,1,1}$.  In the first case, the normalizer of $H$ in $\GL(V)$
is $S_n^2$; in the remaining cases, it is $S_n$.
As representations
of $S_n^2$ or $S_n$ respectively are all self-dual, if $g\in G$ or any scalar multiple thereof
satisfies the Reid-Tai condition, all eigenvalues of $g$ must be $1$ except for a single $-1$.
By the classification of reflection groups, the only one of these possibilities which can
actually occur is the case $V=V_{m-1,1}$ which corresponds to the Weyl group of type $A_{m-1}$.

For case (5), there are only finitely many possibilities for $K$, and for each $K$,
we have $\log |H| \le r \log |\tilde K|$, while the minimal dimension of a faithful irreducible representation of $H$ grows exponentially.  Thus, (\ref{dim.lower.bound}) gives an upper bound on $\dim V$.
The theorem follows from Lemma~\ref{Jordan}.

\end{say}

\section{Quotients of Abelian varieties}

Let us now specialize to  the case when $X=A$ is an Abelian variety
 and $G$ a finite group acting on $A$.
For any $x\in A$, 
the dual of the tangent space $T_xA$ can be canonically identified with
$H^0(A, \Omega_A)$. By Hodge theory
the representation of $\aut_x(A)$ on $H^1(A,\q)\otimes\c$
is isomorphic to the direct sum of the dual
representations on
$H^0(A, \Omega_A)$ and on $H^1(A,\o_A)$.

\begin{defn}
A pair $(G,V)$ is of \emph{AV-type} if $V\oplus V^*$ is isomorphic to 
the complexification of a rational representation of $G$.
\end{defn}

We have the following elementary proposition.

\begin{prop}
\label{few-choices}
Let $(G,V)$ denote a non-RT pair of AV-type.
Let $G_1\subset G$ be a subgroup and $V_1\subset V$ a
$G_1$-subrepresentation such that $(G_1,V_1)$ is 
a basic non-RT pair.
Let $g_1\in G_1$ be exceptional  for $V$, and
let $S_1$ denote the set of eigenvalues of $g_1$ acting on $V_1$,
excluding $1$.   Then every element of $S_1$ is a root of unity whose order lies in
\begin{equation}
\label{possible-orders}
\{2,3,4,5,6,7,8,10,12, 14, 18\}.
\end{equation} 
If $|S_1|>1$, then $S_1$ is one of the following:
\begin{enumerate}
\item $\{e(1/6),e(1/3)\}$.
\item $\{e(1/6),e(1/2)\}$.
\item $\{e(1/6),e(2/3)\}$.
\item $\{e(1/3),e(1/2)\}$.
\item $\{e(1/8),e(3/8)\}$.
\item $\{e(1/8),e(5/8)\}$.
\item A subset of $\{e(1/12),e(1/4),e(5/12)\}$.
\end{enumerate}

\end{prop}

\begin{proof}
Let $\Sigma\subset \Gal(\bar\Q/\Q)$ 
be the set of automorphisms $\sigma$ such that
$V_1^\sigma$ is a $G_1$-subrepresentation of $V$.
As $V\oplus V^*$ is $\Gal(\bar\Q/\Q)$-stable, $\Sigma\cup c\Sigma=\Gal(\bar\Q/\Q)$, where
$c$ denotes complex conjugation.
Let $e(r_1)$ be an element of $S_1$, set
$$S_0=\{\sigma(e(r_1))\mid \sigma\in\Sigma\},$$
and let $r_1,\ldots,r_k$ denote pairwise distinct rational numbers in $(0,1)$ such that 
$S_0 = \{e(r_1),\ldots,e(r_k)\}$.
The $r_i$ have a common denominator $d$, and we write $a_i = d r_i$.
As $\age(g_1)<1$,
$$d > a_1+\cdots+a_k \ge 1+2+\cdots+k\ge \binom{k+1}{2} \ge \binom{\phi(d)/2+1}{2} \ge
\frac{\phi(d)^2}{8}.$$
On the other hand,
$$\phi(d) = d\prod_{p\mid d}\frac{p-1}{p} \ge \frac{d}{3}\prod_{p\mid d, p\ge 5} p^{\frac{\log 4}{\log 5}-1}
\ge \frac{d^{\frac{\log 4}{\log 5}}}{3}.$$
Thus, $d < 372$, and an examination of cases by machine leads to the conclusion
$d$ belongs to the set (\ref{possible-orders}).

If $\alpha$ and $\beta$ are two distinct elements of $S_1$,  then there exists 
$\Sigma$ for which $\Sigma(\alpha)\cup\Sigma(\beta)$ satisfies
the Reid-Tai condition.  On the other hand, $\alpha^f=\beta^f=1$ for some $f\le 126$.
We seek to classify triples of integers $(a,b,f)$, $0<a<b<f\le 126$,
for which there exists a subset $\Sigma\subset\Gal(\bar\Q/\Q)$ with 
$\Sigma\cup c\Sigma = \Gal(\bar\Q/\Q)$ for which
$\Sigma(e(a/f))\cup\Sigma(e(b/f))$ satisfies the Reid-Tai condition.
A machine search for such triples is not difficult and
 reveals that the only possibilities are given by the first six cases of
the proposition together with the three pairs obtained by omitting a single element
from $\{e(1/12),e(1/4),e(5/12)\}$.  The proposition follows.
\end{proof}

\begin{say}[Proof of Theorem~\ref{at.least.four}]
If $A$ is a simple Abelian variety of dimension $\geq4$, $V = T_0 A$, and
$G$ is a finite automorphism group of $A$ which  constitutes an exception to the
statement of the theorem, then $(G,V)$ is a non-RT pair of AV-type.  For every $g\in G$ and every integer $k$, the identity component of the kernel of $g^k-1$, regarded as an endomorphism of $A$, 
is an Abelian subvariety of $A$ and therefore either trivial or equal to the whole of $A$.
It follows that all eigenvalues of $g$ are roots of unity of the same order.

Let $g$ denote an exceptional element.
Let $S_g$ be the multiset of eigenvalues of $g$ on $V$ and $S_g\cup \bar S_g$
the multiset of eigenvalues of $g$ on $V\oplus V^*$.  Then $S_g\cup\bar S_g$ can be
partitioned into a union (in the sense of multisets) of $\Gal(\bar\Q/\Q)$-orbits of roots of unity.
Thus $S_g$ can be partitioned into subsets
$S_{g,i}$ such that for each $X_{g,i}$ either $S_{g,i}\cup \bar S_{g,i}$ or
$S_{g,i}$ itself is a single $\Gal(\bar\Q/\Q)$-orbit.  
If each $S_{g,i}$ is written as a set $\{e(r_{i,1}),\ldots,e(r_{i,j_i})\}$ where the $r_{i,j}$ lie 
in $(0,1)$, then for some $i$, the mean of the values $r_{i,j}$ is less than $\frac14$.  
As every root of unity is Galois-conjugate to its
inverse and $\age(g) < 1$,
if $S_{g,i}$ consists of a single $\Gal(\bar\Q/\Q)$-orbit, then $S_{g,i}$ is $\{1\}$ or $\{-1\}$, in which case
all eigenvalues of $S_g$ are equal, so the mean of the values $r_{i,j}$ is always $\ge \frac12$.
For each $n>2$ in the set (\ref{possible-orders}), we present the set 
$S=\{e(r_1),e(r_2),\ldots,e(r_{\phi(n)/2})\}$ such that $S\cup\bar S$ contains all primitive $n$th roots
of unity and $\sum_j r_j$ is minimal, $r_j\ge 0$:
\goodbreak

\begin{center}
\setlength{\extrarowheight}{-2.75pt}
\renewcommand\arraystretch{1.6}
\begin{tabular}{|l|c|c|c|}
\hline 
$n$&\emph{$\phi(n)/2$}&\emph{Values of $r_j$}&\emph{Mean of $r_j$} \\ 
\hline
$3$&$1$&$\frac13$&$\frac13$ \\
$4$&$1$&$\frac14$&$\frac14$ \\
$5$&$2$&$\frac15,\frac25$&$\frac3{10}$ \\
$6$&$1$&$\frac16$&$\frac16$ \\
$7$&$3$&$\frac17,\frac27,\frac37$&$\frac27$ \\
$8$&$2$&$\frac18,\frac38$&$\frac14$ \\
$9$&$3$&$\frac19,\frac29,\frac49$&$\frac7{27}$ \\
$10$&$2$&$\frac1{10},\frac3{10}$&$\frac15$ \\
$12$&$2$&$\frac1{12},\frac5{12}$&$\frac14$ \\
$14$&$3$&$\frac1{14},\frac3{14},\frac5{14}$&$\frac3{14}$ \\
$18$&$3$&$\frac1{18},\frac5{18},\frac7{18}$&$\frac{13}{54}$ \\
\hline
\end{tabular}
\renewcommand\arraystretch{1.0}
\smallskip \\ Table 1
\end{center}

Inspection of this table reveals that the mean of the values $r_j$ is less than $\frac 14$, 
only if $n=18$, $n=14$, $n=10$, or $n=6$.   In the first two cases, the condition $\dim A\ge 4$
implies there must be at least two subsets $S_{g,i}$ in the partition, which implies $\dim A\ge 6$.
As the mean of the $r_j$ exceeds $\frac16$ for $n=14$ and $n=18$, this is impossible.
If $n=6$, all the eigenvalues of $g$ must be $\frac 16$, and there could be as many as five.
However, an abelian variety $A$ with an automorphism which acts as the scalar $e(1/6)$
on $T_0 A$ is of the form $\C^g/\Lambda$, where $\Lambda$ is a torsion-free 
$\Z[e(1/6)] = \Z[e(1/3)]$-module
with the inclusion $\Lambda\to \C^g$ equivariant with respect to $\Z[e(1/3)]$.  Every
finitely generated torsion-free module over $\Z[e(1/3)]$ is free (as $\Z[e(1/3)]$ is a PID), 
so $A$ decomposes as a product of elliptic curves with CM by $\Z[e(1/3)]$, contrary to hypothesis.
If $n=10$, the only possibility is that $\dim A=4$, and the eigenvalues of $g$
are $e(1/10),e(1/10),e(3/10),e(3/10)$.  Again, $\Z[e(1/5)]$ is a PID, so $A = \C^4/\Lambda$,
where $\Lambda\cong \Z[e(1/5)]\oplus \Z[e(1/5)]$.  Let $\Lambda_1\subset \Lambda$ denote the 
first summand.  As $A$ is simple, the $\C$-span of
$\Lambda_1$ must have dimension $>2$.  However, letting $\lambda\in\Lambda_1$ be a generator,
we can write $\lambda$ as a sum of two eigenvectors for $e(1/10)\in \Z[e(1/5)]$.
Every element of $\Z[e(1/5)]\lambda$ is then a complex linear combination of these $2$ eigenvectors.
\end{say}


\begin{lemma}
\label{imprimitive.ok}
Let $(G,V)$ denote a non-RT pair of AV-type.
Let $G_1\subset G$ be a subgroup and $V_1\subset V$ a
$G_1$-subrepresentation such that $(G_1,V_1)$ is 
an imprimitive basic non-RT pair.  Then $(G_1,V_1)$ is a complex reflection group.

\end{lemma}

\begin{proof}
By Proposition~\ref{Prod}, $V_1$ decomposes as a direct sum of lines which $G_1$ permutes.
Let $g_1$ be an element of $G_1$ which is exceptional for
$V_1$.  By Proposition~\ref{Prod}, after renumbering the $L_i$, 
$g_1$ interchanges $L_1$ and $L_2$ and
stabilizes all the other $L_i$.
The eigenvalues of $g_1$ acting on $L_1\oplus L_2$ are therefore
of the form $e(r)$ and $e(r+1/2)$ for some $r\in [0,1/2)$.  By Proposition~\ref{few-choices},
this means $r=0$, $r=1/6$, or $r=1/8$.  In the first case, $g$ might have an additional
eigenvalue $e(1/6)$ or $e(1/3)$ on one of the lines $L_i$, $i\ge 3$ and fix all the remaining
lines pointwise.  In all other cases, $g_1$ must fix $L_i$ pointwise for $i\ge 3$.
If $g_1$ has eigenvalues $-1,1,e(1/6),1,\ldots,1$, eigenvalues $-1,1,e(1/3),1\ldots,1$, 
eigenvalues $e(1/6),e(2/3),1,\ldots,1$, or eigenvalues $e(1/8),e(5/8),1,\ldots,1$, then
$g_1^2$ is again exceptional but stabilizes all of the lines $L_i$, which is impossible
by Proposition~\ref{Prod}.
In the remaining case, $g_1$ has eigenvalues $-1,1\ldots,1$, so $g_1$ is a reflection,
and $G_1$ is a complex reflection group.
\end{proof}

\begin{thm}
\label{AV.case}
Let $(G,V)$ be a non-RT pair of AV-type,
$G_1\subset G$ a subgroup, and $V_1\subset V$ a
$G_1$-subrepresentation such that $(G_1,V_1)$ is 
a basic non-RT pair.   If $(G_1,V_1)$ is not of reflection type, then $\dim V_1 = 4$,
and $G_1$ is contained in the reflection group $G_{31}$ in the Shephard-Todd classification.
\end{thm}

\begin{proof}
Let $g_1\in G_1$ be an exceptional element for $V$.
If $g_1$ or any of its powers is a pseudoreflection on $V_1$, then
$G_1$ is generated by pseudoreflections and is therefore a reflection group on $V_1$.  
We may therefore assume that every power of $g_1$ which is non-trivial on
$V_1$ has at least two non-trivial eigenvalues in its action on $V_1$.
Also, by Lemma~\ref{imprimitive.ok}, $(G_1,V_1)$ may be assumed primitive.

We consider first the case that $S_1$ consists of a single element of order $n$.
By a well-known theorem of Blichfeldt (see, e.g., \cite[5.1]{Cohen})
a non-scalar element in a primitive group cannot have all of its eigenvalues contained in an arc of length $\pi/3$.  If follows that $n\le 5$.  If $n=5$,
$V$ 
contains at least two different fifth roots of unity, and as $g_1$ is exceptional on $V$,
it follows that the multiplicity of the non-trivial eigenvalue of $g_1$ on $V_1$ is $1$, contrary to
hypothesis.  If $n=4$, then by \cite{wales}, the eigenvalues $1$ and $i$
have the same multiplicity (which must be at least $2$), 
and by \cite{korlyukov}, $\dim V_1$ is a power of $2$.  As the multiplicity of $i$ is at most $3$, the only
possibility is that $\dim V_1 = 4$ and the eigenvalues of $g_1$ are $1,1,i,i$.  The classification
of primitive $4$-dimensional groups \cite{blichfeldt} shows that the only such groups containing
such an element are contained in the group $G_{31}$ in the Shephard-Todd classification.
If $n=3$, the multiplicity of the eigenvalue $e(1/3)$ must be $2$, and
\cite{wales} show that there are only two possible examples, one in dimension $3$
(which is projectively equivalent to the Hessian reflection group $G_{25}$)
and one in dimension $5$ (which is projectively equivalent to the reflection group $G_{33}$).
The case $n=2$ does not arise, since the non-trivial eigenvalue multiplicity is $\ge 2$.

Thus, we need only consider the cases that $|S_1|\ge 2$.  The possibilities for the
multiset of non-trivial eigenvalues of $g_1$ acting on $V_1$ which are consistent
with $g_1$ being exceptional for an AV-pair $(G_1,V)$ are as follows:

\renewcommand{\theenumi}{\alph{enumi}}
\begin{enumerate}
\item $e(1/6),e(1/3)$
\item $e(1/6),e(1/6),e(1/3)$
\item $e(1/6),e(1/6),e(1/6),e(1/3)$
\item $e(1/6),e(1/3),e(1/3)$
\item $e(1/6),e(1/2)$
\item $e(1/6),e(1/6),e(1/2)$
\item $e(1/6),e(2/3)$
\item $e(1/3),e(1/2)$
\item $e(1/8),e(3/8)$
\item $e(1/8),e(5/8)$
\item $e(1/12),e(1/4)$
\item $e(1/12),e(5/12)$
\item $e(1/4),e(5/12)$
\item $e(1/12),e(1/4),e(5/12)$
\end{enumerate}
\renewcommand{\theenumi}{\arabic{enumi}}

Cases (a), (d), (e), (g), (h), (k), and (m) are ruled out because no power of $g_1$
may be a pseudoreflection.  In case (f), $g_1^2$ has two non-trivial eigenvalues, both equal to $e(1/3)$,
and we have already treated this case.  Likewise, cases (j) and (l) are subsumed in our analysis of
the case that there are two non-trivial eigenvalues, both equal to $i$.

For the four remaining cases, we observe that the conjugacy class of $g_1$ generates
the non-abelian group $G_1$, so $g_1$ fails to commute with some conjugate $h_1$.
The group generated by $g_1$ and $h_1$ fixes a subspace $W_1$ of $V_1$ of codimension 
at most $6$, $8$, $4$, and $6$ in cases (b), (c), (i), and (n) respectively which 
$g_1$ and $h_1$ fix pointwise.  Let $U_1\subset V_1/W_1$ denote a space on which
$\langle g_1,h_1\rangle$ acts irreducibly and on which $g_1$ and $h_1$ do not commute.
The non-trivial eigenvalues of $g_1$ and $h_1$ on $U_1$ form subsets of the non-trivial
eigenvalues of $g_1$ and $h_1$ on $V_1$, and the action of $\langle g_1,h_1\rangle$
on $U_1$ is primitive because the eigenvalues of $g_1$ do not include a coset of
any non-trivial subgroup of $\C^\times$.

We claim that if $\dim U_1>1$, all the non-trivial eigenvalues of $g_1$ on $V_1$
occur already in $U_1$.   In cases (i) and (n), we have already seen that no proper subset of
indicated sets of eigenvalues can appear, together with the eigenvalue $1$ with some multiplicity,
in any primitive irreducible representation.  In cases (b) and (c), Blichfeldt's $\pi/3$ theorem
implies that if the eigenvalues of some element in a primitive representation of a finite group
are $1$ with some multiplicity, $e(1/6)$ with some multiplicity,
and possibly $e(1/3)$, then $e(1/3)$ must actually appear.  Therefore, the factor $U_1$
must have $e(1/3)$ as eigenvalue, and every other irreducible factor of $V_1$ must be $1$-dimensional.  If no eigenvalue $e(1/6)$ appears in $g_1$ acting on $U_1$,
then $g_1^3$ and $h_1^3$ commute.  If all conjugates of $g_1^3$ commute, then $G$ has a
normal abelian subgroup.  Such a subgroup must consist of scalar elements of $\End(V_1)$, but this is
not possible given that at least one eigenvalue of $g_1^3$ on $V_1$ is $1$ and at least one eigenvalue is $-1$.  Without loss of generality, therefore, we may assume that $g_1^3$ and $h_1^3$ fail to commute.  It follows that both $e(1/3)$ and $e(1/6)$ are eigenvalues of $g_1$ on $U_1$.
As case (a) has already been disposed of, the multiplicity of $e(1/6)$ as an eigenvalue of $g_1$ on $U_1$ is at least $2$.  In the case (b), this proves the claim.  Once it has been shown that there are
no solutions of type (b), it will follow that the eigenvalue $e(1/6)$ must appear with multiplicity $3$,
which proves the claim for (c).  

Finally, we show that for each of the cases (b), (c), (i), and (n) there is no finite group $G_1$
with a primitive representation $U_1$ and an element $g_1$ 
whose multiset of non-trivial eigenvalues is
as specified.   First we consider whether $G_1$ can stabilize a non-trivial tensor decomposition of
$U_1$.  The only possibilities for $g_1$ respecting such a decomposition are case (b)
with eigenvalues $1,e(1/6),e(1/6),e(1/3)$ decomposing as a tensor product of two representations
with eigenvalues $1,e(1/6)$ and case (c) with eigenvalues $1,1,e(1/6),e(1/6),e(1/6),e(1/3)$
decomposing as a tensor product of representations with eigenvalues $1,e(1/6)$ and
$1,1,e(1/6)$.  As $U_1$ is a primitive representation of $G_1$, Blichfeldt's theorem rules out both possibilities.  

Next we rule out the possibility that $G_1$ normalizes a tensor decomposition with $g_1$ permuting tensor factors non-trivially.  Given that $\dim U_1\le 8$, this could only happen
if there are two or three tensor factors, each of dimension $2$.   It is easy to see that if
$T_1,\ldots,T_n$ are linear transformations on a vector space $V$, the transformation
on $V^{\otimes n}$ defined by $v_1\otimes\cdots\otimes v_n\mapsto T_n(v_n)\otimes T_1(v_1)\otimes\cdots\otimes T_{n-1}(v_{n-1})$ has the same trace as $T_1 T_2\cdots T_n$.
It follows that any unitary transformation $T$ 
on $V^{\otimes n}$ which normalizes the tensor decomposition
but permutes the factors nontrivially satisfies
$$|\tr(T)| \le (\dim V)^{n-1},$$
with equality only if the permutation is a transposition $(ij)$, $T_i T_j$ is scalar, and all other
factors $T_i$ are scalar; in particular, equality implies that $T^2$ is scalar.
The following table gives for each case 
the absolute value of the trace of $g_1$ acting on $U_1$ in terms
of the dimension of $U_1$:
\goodbreak
\begin{center}
\begin{tabular}{|l|l|l|l|l|l|l|}
\hline 
\emph{Case}&$3$&$4$&$5$&$6$&$7$&$8$ \\ 
\hline
(b)&$\sqrt3$&$2$        &$\sqrt7$&$2\sqrt3$   &                   & \\
(c)&               &$\sqrt7$&$3$        &$\sqrt{13}$&$\sqrt{19}$&$3\sqrt3$ \\
(i)&$\sqrt3$ &$\sqrt6$&               &                    &                    & \\
(n)&$2$        &$\sqrt5$&$2\sqrt2$&$\sqrt{13}$&                  & \\
\hline
\end{tabular}
\smallskip \\ Table 2
\end{center}
In each case, except (b) and $\dim U_1=4$, $\tr(T)$ violates the inequality, and in this case,
$T^2$ is not scalar.

Let $H_1$ denote a characteristically simple normal subgroup of $G_1$.  As $G_1$ does not normalize a tensor decomposition, $U_1$ is an irreducible representation of $H_1$.  
Either $H_1$ is the product of an extraspecial $p$-group $H_p$ and a group $Z$ of scalars or
$H_1$ is a central extension of a product of mutually isomorphic finite simple groups by a scalar group $Z$.  As $\dim U_1\le 8$, in the former case, $|H_p|\in\{2^3,2^5,2^7,3^3,5^3,7^3\}$.
In the latter case, either $\bar H_1 = H_1/Z$ is isomorphic to $K^r$ for some finite simple group $K$,
and $r=1$ since $G_1$ does not normalize a tensor decomposition.
For a list of possibilities for $H_1$, we use the tables of Hiss and Malle \cite{hiss-malle}, which are
based on the classification of finite simple groups.  Note that primitive groups were classified up through dimension $10$ before the classification of finite simple groups
was available (see, e.g., \cite{Feit70}, \cite{Feit75} and the references therein).
The following table enumerates the possibilities for $\bar H_1$, where representation numbering
is that of \cite{Atlas} and asterisks indicate a Stone-von Neumann representation:

\begin{center}
\begin{tabular}{|l|c|c|c|c|c|c|c|}
\hline
&\multicolumn{7}{c|}{Representation Degree} \\
\hline 
Group&2&3&4    &5&6           &7&8 \\
\hline
$(\Z/2\Z)^2$&$*$&&&&&& \\
$(\Z/3\Z)^2$&&$*$&&&&& \\
$(\Z/2\Z)^4$&&&$*$&&&& \\
$(\Z/5\Z)^2$&&&&$*$&&& \\
$(\Z/7\Z)^2$&&&&&&$*$& \\
$A_5$    &6&3&4,\,8&5&9            &  &   \\
$(\Z/2\Z)^6$&&&&&&&$*$ \\
$L_2(7)$&  &2&7    &  &4,\,9        &5&6,\,11 \\
$A_6$     & &14&8  &2&16,\,19   & &4,\,10  \\
$L_2(8)$& &     &     &  &            &2,\,3&6 \\
$L_2(11)$&&   &      &2&9          &   & \\
$L_2(13)$&&   &      & &10        & 2& \\
$L_2(17)$&&   &      & &             &   &12 \\
$A_7$     & &    &10  & &2,\,17,\,24& & \\
$U_3(3)$&&    &       & &2          &3,\,4& \\
$A_8$    &  &    &       & &            &2&15 \\
$L_3(4)$& &    &       & &41       &  &19 \\
$U_4(2)$& &   &21  &2&4        &   & \\
$A_9$   &  &    &       & &            &  &2,\,19 \\
$J_2$   &   &    &       & &22      &   & \\
$S_6(2)$& &   &       & &            &2&31 \\
$U_4(3)$& &   &       & &72      &  & \\
$O_8^+(2)$&&&&&&&54 \\
\hline
\end{tabular}
\smallskip \\ Table 3
\end{center}

For the finite simple groups $\bar H_1$, we consult character tables \cite{Atlas}.  This is easy to do so, since only a few of the characters in Table 3 take values whose absolute values are large
enough to appear in Table 2.  There are no cases where an element of order $6$ has a character
absolute value as given in row (b) or (c) of Table 2; an element of order $8$ has an absolute value as given by row (i); or an element of order $12$ has an absolute value as given by row (n).

For the case that $H_1$ is an extraspecial  $p$-group, every non-zero character value is
an integral power of $\sqrt{p}$.  This was proved for $p>2$ by Howe \cite[Prop.~2(ii)]{Howe}.
For lack of a reference for $p=2$, we sketch a proof which works in general.
The embedding $H_1\to\GL(U_1)$ is a Stone-von Neumann representation with central character
$\chi$.  Let $GG_1$ denote the group of pairs $(g_1,g_2)\in G_1^2$ such that 
$g_1 H_1 = g_2 H_1$.  There is a natural action of $GG_1$ on $\C[H_1]$ given by
$$(g_1,g_2)([h_1]) = [g_1 h_1 g_2^{-1}].$$
The restriction of this representation to $H_1^2 \subset GG_1$ is of the form 
$\bigoplus V_i\boxtimes V_i^*$, where the sum is taken over all irreducible representations
$V_i$ of $H_1$.  The factor $U_1\boxtimes U_1^*$
is the $\chi\boxtimes \chi^*$ eigenspace of the center $Z^2$ of $H_1^2$ acting on $\C[H_1]$,
where $\chi$ is the central character of $Z$ on $U_1$.  As the action of $GG_1$ on this
eigenspace of $\C[H_1]$ extends  the irreducible representation of $H_1^2$ on
$U_1\boxtimes U_1^*$, any other extension of $(H_1^2,U_1\boxtimes U_1^*)$ to
$GG_1$ is projectively equivalent to this one.  The particular extension we have in mind is
obtained by letting $(g_1,g_2)\in GG_1$ act on $U_1\boxtimes U_1^*$ according to the action of
$g_1$ on $U_1$ and the action of $g_2$ on $U_1^*$ coming from the inclusion $G_1\subset
\GL(U_1)$.  From this it is easy to see that the character value of $(g_1,g_1)$ on
each $Z^2$-eigenspace of $\C[H_1]$ is either $0$ or $|\bar H_1^{g_1}|$.  As
$\bar H_1^{g_1}$ is a vector space over $\Z/p\Z$, $\tr(g_1|U_1)\tr(g_1|U_1^*)$
is either $0$ or an integer power of $p$.  Consulting Table 2, we see that this rules out every possibility except a character value $2$ and $\dim U_1 = 4$.  This can actually occur, but not with the eigenvalues of case (b).

\end{proof}

\section{Examples}

We conclude with some examples to illustrate various aspects of the classification given above.
We begin with some examples from group theory.

In principle all non-RT pairs can be built up from basic pairs, by reversing the operations
which led to constructing basic pairs in the first place, i.e. by replacing $G$ by an extension
$\tilde G$ of $G$ whose image in $\Aut(V)$ is the same as that of $G$; by combining
$(G,V_1)$ and $(G,V_2)$ to give the pair $(G,V_1\oplus V_2)$ (which may or may not be non-RT); 
and by replacing
$(G,V)$ by $(G',V)$, where $G'$ lies between $G$ and its normalizer in $\Aut(V)$.
To illustrate that, we observe that all non-RT pairs of the form
$((\Z/2\Z)^n\rtimes H,\C^n)$, where $H\subset S_n$ is a transitive group, arise from the basic
non-RT pair $(\Z/2\Z,\C)$.  This accounts for the series of Weyl groups of type $B_n/C_n$,
but not for the Weyl groups of type $D_n$, which are primitive.
This construction can be used more generally to build non-RT pairs of the form 
$(G^n\rtimes H, V^n)$ starting with a non-RT pair $(G,V)$ and a transitive
permutation group $H\subset S_n$.  

It may happen that a basic non-RT pair $(G,V)$ of reflection type nevertheless fails to have an exceptional element which is a scalar multiple of a pseudoreflection.  Consider the case 
$G = U_4(2)\times \Z/3\Z$ and $V$ is a faithful irreducible $5$-dimensional representation of $G$.  Then $G$ has an exceptional element $g$ whose eigenvalues are $1,1,1,e(1/3),e(1/3)$ and whose conjugacy class generates $G$.  It has another element $h$ with eigenvalues $1,-1,-1,-1,-1$ which is not exceptional.  As $-h$ is a reflection, it is easy to see that $U_4(2)\times \Z/2\Z$ is a $5$-dimensional reflection group (in fact, it is $G_{33}$ in the Shephard-Todd classification), and of course this reflection group is projectively equivalent to $(G,V)$.

There really does exist a primitive $4$-dimensional non-RT pair $(G,V)$ which is not of reflection type.
By \cite{blichfeldt}, there is a short exact sequence
$$0\to I_4\to G_{31}\to S_6\to 0$$
where $I_4$ is the central product of $\Z/4\Z$ and any extraspecial $2$-group of order $32$.
The group $S_6$ contains two non-conjugate subgroups isomorphic to $S_5$, whose inverse
images in $G_{31}$ are primitive.  One is the reflection group $G_{29}$, and one contains
elements with eigenvalues $1,1,i,i$.  The question arises as to whether these
two groups are conjugate.  
The character table of $G_{29}$, provided by the software package CHEVIE
\cite{GHLMP}, reveals that this group has two faithful $4$-dimensional representations.
One has reflections and the other has elements with spectrum $1,1,i,i$.
It follows that $G_{29}$ with respect to this non-reflection representation, or equivalently,
the non-reflection index-$6$ subgroup of $G_{31}$, gives the desired example.
This example (in fact all of $G_{31}$) can actually be realized inside $\GL_4(\Z[i])$,
as shown in \cite{blichfeldt}.

The set of projective equivalence classes of basic non-RT pairs which 
are of AV-type is infinite, as is the set of basic non-RT pairs which are not.
We have already mentioned the Weyl groups of type $D_n$ as examples of the first kind;
the reflection groups $(\Z/k\Z)^{n-1}\times S_n$ are never of AV-type if $k > 4$.

We conclude with some geometric examples.

If $(G,V)$ is a non-RT pair and $V = V_0\otimes_{\Q}\C$ for some rational representation
$V_0$ of $G$, then there exists an Abelian variety $A$, and a homomorphism
$G\to \Aut(A)$ such that $A/G$ is uniruled and the Lie algebra of $A$ is isomorphic to $V$
as $G$-module.  Indeed, we may choose any integral lattice $\Lambda_0\subset V_0$ which
is $G$-stable and define $A = \Hom(\Lambda_0,E)$ for any elliptic curve $E$.  If $V$ is irreducible,
then $A/G$ is rationally connected.  This includes all examples  where $G$ is the Weyl group
of a root system and $V_0$ is the $\Q$-span of the root system.  When $\Lambda_0$
is taken to be the root system, the quotients $A/G$ are in fact weighted projective spaces by
a theorem of E.~Looijenga \cite{looijenga} which was one of the motivations for this paper.

Let $\Lambda_0$ denote the ($12$-dimensional) Coxeter-Todd lattice, which we regard
as a free module of rank $6$ over $R:=\Z[e(1/3)]$.  Let $G = G_{34}$ denote the group of
$R$-linear isometries of this lattice.  If $E$
denotes the elliptic curve over $\C$ with complex multiplication by $R$, then
$\Hom_R(\Lambda_0,E)/G$ is rationally connected.  The group $G$ is a reflection group
but not a Weyl group,
and we do not know whether this variety is rational or even unirational.

We have already observed that there is a $4$-dimension basic non-RT pair $(G,V)$
which is not of reflection type and such that $G \subset \GL_4(\Z[i])\subset \GL(V)$.  
If $E$ denotes the elliptic curve with CM by $\Z[i]$, $E^4/G$ is
rationally connected.  Again, we do not know about rationality or unirationality.

Let $A$ be an abelian variety with complex multiplication by $\Z[e(1/7)]$ with CM type
chosen so that some automorphism $g$ of order $7$ has eigenvalues $e(1/7),e(2/7),e(3/7)$
acting on $T_0 A$.  Then $A/\langle g\rangle$ is rationally connected, but once again we do not know whether it is rational or unirational.

Let $E$ be any elliptic curve, $A=E^3$, and $G = S_3\times \{\pm1\}$.  Let the factors $S_3$ and $\{\pm1\}$ of $G$ act on $A$ by permuting factors and by multiplication respectively.
If $V=T_0 A = \C^3$ and $W$ denotes the plane in which coordinates sum to zero, then
the images of $G$ in $W$ and in $V/W$ are reflection groups.
Thus $A^{RT}_{stab} = A$, so $A/G$ is rationally connected by Corollary~\ref{filtered.rc}.  
Yet again, we do not know about the rationality or unirationality of the quotient.

There are imprimitive basic non-RT triples of arbitrarily large degree which can be realized
by automorphism groups of Calabi-Yau varieties but not as automorphism groups of
Abelian varieties.  For example,
$$G_{n+2,1,n} = (\Z/(n+2)\Z)^n\rtimes S_n$$
acts on the $n$-dimensional Fermat hypersurface $x_0^{n+2}+\cdots+x_{n+1}^{n+2}=0$
fixing the coordinates $x_0$ and $x_1$ and therefore the point 
$$P = (1:e^{\frac{\pi i}{n+2}}:0:\cdots:0).$$
The action of $G$ on the tangent space to $P$ gives the reflection representation of
$G_{n+2,1,n}$. 

\bibliography{refs3}

\end{document}